\documentclass[10pt, a4paper]{amsart}

\usepackage{amsmath}
\usepackage{amsfonts, mathtools}
\usepackage{amsthm, upref, url}
\usepackage{graphicx, enumitem}
\usepackage[usenames, dvipsnames]{color}
\usepackage[colorlinks=true,
        raiselinks=true,
        pagecolor=magenta,
        linkcolor=MidnightBlue,
        citecolor=Mahogany,
        urlcolor=PineGreen,
        pdfauthor={Debasish Chatterjee},
        pdftitle={An Excursion-Theoretic Approach to Stability of Discrete-Time Stochastic Hybrid Systems},
        pdfkeywords={Markov processes, stochastic stability},
        pdfsubject={Technical Report},
        plainpages=false]{hyperref}
\usepackage[bitstream-charter]{mathdesign}
\DeclareSymbolFont{usualmathcal}{OMS}{cmsy}{m}{n}
\DeclareSymbolFontAlphabet{\mathcal}{usualmathcal}

\DeclareMathOperator*{\esssup}{ess\,sup}

\newcommand{\rr}{\mathbb{R}}
\newcommand{\mcal}{\mathcal}
\newcommand{\id}{\mathrm{id}}
\newcommand{\iss}{\textsc{iss}}
\newcommand{\ifs}{\textsc{ifs}}
\newcommand{\msf}{\mathsf}
\newcommand{\tauk}{\ensuremath{\tau_{K}^{\vphantom{T}}}}
\newcommand{\etalchar}[1]{$^{#1}$}
\newcommand{\sts}{\mathcal{S}}
\newcommand{\dup}{{\mathrm{d}}}
\newcommand{\eq}{\begin{equation}}
\newcommand{\en}{\end{equation}}
\newcommand{\stp}{\mathbb{T}}
\newcommand{\Let}{\coloneqq}

\newcommand{\R}{\ensuremath{\mathbb{R}}}
\newcommand{\N}{\ensuremath{\mathbb{N}}}
\newcommand{\Nz}{\ensuremath{\mathbb{N}_0}}
\newcommand{\posR}{\ensuremath{\R_{\ge 0}}}
\newcommand{\Pm}{\ensuremath{\mathsf P}}
\newcommand{\EE}{\ensuremath{\mathsf E}}
\newcommand{\transp}{\ensuremath{^{\scriptscriptstyle{\mathrm T}}}}
\newcommand{\sigalg}{\ensuremath{\mathfrak{F}}}
\newcommand{\mx}{\ensuremath{\vee}}
\newcommand{\mn}{\ensuremath{\wedge}}
\newcommand{\lra}{\ensuremath{\longrightarrow}}
\newcommand{\xz}{\ensuremath{x_0}}
\newcommand{\therex}{\ensuremath{\exists\,}}
\newcommand{\setmin}{\ensuremath{\setminus}}
\newcommand{\cardP}{\ensuremath{\mathrm{N}}}
\newcommand{\PSet}{\ensuremath{\mathcal{P}}}
\newcommand{\eps}{\ensuremath{\varepsilon}}
\newcommand{\fa}{\ensuremath{\forall\,}}
\renewcommand{\le}{\ensuremath{\leqslant}}
\renewcommand{\ge}{\ensuremath{\geqslant}}
\newcommand{\ClassK}{\ensuremath{\mathcal{K}}}
\newcommand{\ClassKinfty}{\ensuremath{\mathcal{K}_{\infty}}}
\newcommand{\ClassKL}{\ensuremath{\mathcal{KL}}}

\newcommand{\Ex}[1]{\ensuremath{\mathsf{E}\!\left[\vphantom{\big|}#1\vphantom{\big|}\right]}}

\newcommand{\indic}[1]{\ensuremath{\boldsymbol{1}_{#1}}}
\newcommand{\mrm}[1]{\ensuremath{\mathrm{#1}}}

\newcommand{\Lp}[1]{\ensuremath{\boldsymbol{L}_{#1}}}
\newcommand{\secref}[1]{\S\ref{#1}}
\newcommand{\epower}[1]{\ensuremath{\mathrm{e}^{#1}}}
\newcommand{\norm}[1]{\ensuremath{\left\lVert #1 \right\rVert}}

\newcommand{\RemarkEnd}{\hspace{\stretch{1}}{$\vartriangleleft$}}

\newcommand{\DefEnd}{\hspace{\stretch{1}}{$\Diamond$}}
\newcommand{\AssumptionEnd}{\hspace{\stretch{1}}{$\diamondsuit$}}

\allowdisplaybreaks

\numberwithin{equation}{section}
\swapnumbers
\newtheoremstyle{dcstyle}{8pt}{8pt}{\slshape}{}{\bfseries}{.}{ }{}

\theoremstyle{dcstyle}
\newtheorem{theorem}[equation]{Theorem}
\newtheorem{corollary}[equation]{Corollary}
\newtheorem{lemma}[equation]{Lemma}
\newtheorem{proposition}[equation]{Proposition}

\theoremstyle{definition}
\newtheorem{defn}[equation]{Definition}

\theoremstyle{remark}
\newtheorem{remark}[equation]{Remark}

\newtheorem{assumption}[equation]{Assumption}

\makeatletter
\def\tagform@#1{\maketag@@@{\ignorespaces#1\unskip\@@italiccorr}}
\makeatother

\title[Excursion-Theoretic view of Stochastic Hybrid System Stability]{An Excursion-Theoretic Approach to Stability of Discrete-Time Stochastic Hybrid Systems}

\author{Debasish Chatterjee}
\address{ETL I19\\Physikstrasse 3\\ETH Z\"urich\\8092 Z\"urich\\Switzerland\\Ph.~+41-44-632-2326\\Fax:~+41-44-632-1211}
\email{chatterjee@control.ee.ethz.ch}

\author{Soumik Pal}
\address{C-547 Padelford Hall\\Department of Mathematics\\University of Washington, Seattle\\WA~98195\\Ph.~+1-206-543-7832}
\email{soumik@math.washington.edu}



\thanks{Debasish Chatterjee's research is partially supported by the Swiss National Science foundation grant 200021-122072.}

\date{\today}

\begin{document}

	\begin{abstract}
		We address stability of a class of Markovian discrete-time stochastic hybrid systems. This class of systems is characterized by the state-space of the system being partitioned into a safe or target set and its exterior, and the dynamics of the system being different in each domain. We give conditions for $\Lp 1$-boundedness of Lyapunov functions based on certain negative drift conditions outside the target set, together with some more minor assumptions. We then apply our results to a wide class of randomly switched systems (or iterated function systems), for which we give conditions for global asymptotic stability almost surely and in $\Lp 1$. The systems need not be time-homogeneous, and our results apply to certain systems for which functional-analytic or martingale-based estimates are difficult or impossible to get.
	\end{abstract}

	\maketitle

\section{Introduction}
	Increasing complexity of engineering systems in the modern world has led to the hybrid systems paradigm in systems and control theory~\cite{vanderschaftHbk, liberzonbk}. A hybrid system consists of a number of domains in the state-space and a dynamical law corresponding to each domain; thus, at any instant of time the dynamics of the system depends on the domain that its state is in. One would then restrict attention to behavior of the system in individual domains, which is typically a simpler problem. However, understanding how the dynamics in the individual domains interact among each other is necessary in order to ensure smooth operation of the overall system. This article is a step towards understanding the behavior of (possibly non-Markovian) stochastic hybrid systems which undergo excursions into different domains infinitely often. Here we consider the simplest and perhaps the most important hybrid system, consisting of a compact target or safe set and its exterior, with different dynamics inside and outside the safe set. Our objective is to introduce a new method of analysis of systems that are outside the safe set infinitely often in course of their evolution. The analysis carried out here provides a basis for controller synthesis of systems with control inputs---it gives clear indications about the type of controllers to be designed in order to ensure certain natural and basic stability properties in closed loop.

	Let us look at two interesting and practically important examples of hybrid systems with two domains---a compact safe set and its exterior, with different dynamics in each. The first concerns optimal control of a Markov process with state constraints. Markov control processes have been extensively studied; we refer the reader to the excellent monographs and surveys~\cite{bertsekasshreve78, borkarTopicsControlledMC, hernandez-lerma1, hernandez-lerma2} for further information, applications and references. For our purposes here, consider the canonical example of a linear controlled system perturbed by additive Gaussian noise and having probabilistic constraints on the states. A hybrid structure of the controlled system naturally presents itself in the following fashion. Except in the most trivial of cases, computing the constrained optimal control over an infinite horizon is impossible, and one resorts to a a rolling-horizon controller. (Rolling-horizon controllers are considerably popular, for basic definitions, comparisons and references see e.g., \cite{maciejowskibk} in the deterministic context, and~\cite{chl09} and the references therein in the stochastic context.) Computational overheads restrict the size of the window in the rolling-horizon controller, and determine the maximal (typically bounded) region---called the safe set---in which this controller can be active. No matter how good the resulting controller is, the additive nature of the Gaussian noise ensures that the states are subjected to excursions away from the safe set infinitely often almost surely. Once outside the safe set, the rolling-horizon controller is switched off and a recovery strategy is activated, whose task is to bring the states back to the safe set quickly and efficiently. This problem is of great practical interest and a subject of current research, see e.g., \cite{chl09} and the references in them for possible strategies inside the safe set, and~\cite{cccl08} for one possible recovery strategy. Evidently, stability of this hybrid system depends largely on the recovery strategy, since as long as the states stay inside the safe set, they are bounded. However, traditional methods of stability analysis do not work well precisely because of the unlimited number of excursions. Theorem~\ref{maintheorem} of this article addresses this issue, and provides a method of ensuring strong boundedness and stability properties of the hybrid system. Intuitively it says that under the recovery strategy there exists a well-behaved supermartingale until the states hit the safe set, then the system state is bounded in expectation uniformly over time. A complete picture of stability and ergodic properties of a general controlled hybrid system is beyond the scope of the present article, and will be reported elsewhere. We refer the reader to~\cite[Chapter~3]{costaMJLS} for earlier work pertaining to stability of a class of hybrid systems, and to~\cite{meynbk} for stability of general discrete-time Markov processes.

	The second example is one that we shall pursue further in this article, namely, a class of discrete-time Markov processes called iterated function systems~\cite{barnsley88, lasota94} (\ifs{}). They are widely applied, for instance, in the construction of fractals~\cite{lasota94}, in studies on the process of generation of red blood corpuscles~\cite{lasotamyjak02, lasotaszarek04}, in statistical physics~\cite{kifer86}, and simulation of important stochastic processes~\cite{werner05}. Of late they are being employed in key problems of physical chemistry and computational biology, namely, the behavior of the chemical master equation~\cite[Chapter~6]{wilkinson06} (CME), which governs the continuous-time stochastic (Markovian) reaction-kinetics at very low concentrations (of the order of tens of molecules). Invariant distributions, certain finite-time properties, and robustness properties with respect to disturbances of the underlying Markov process are of interest in modeling and analysis of unicellular organisms. It is well-known that the CME is analytically intractable (see~\cite{huisinga07, kurtz08} for special cases), but the invariant distribution of the Markov process can be recovered from simulation of the embedded Markov chain in a computationally efficient way~\cite{miliasforthcoming}. This embedded chain is an \ifs{} taking values in a nonnegative integer lattice. From a biological perspective, good health of a cell corresponds to the \ifs{} evolving in a safe region on an average, despite moderate disturbances to the numbers of molecules involved in the key reactions. However, in most cases compact invariant sets do not exist. It is therefore of interest to find conditions under which, even though there are excursions of the states away from a safe set infinitely often, the \ifs{} is stochastically bounded, or some strong stability properties hold. Theorem~\ref{maintheorem} of this article leads to results (in \secref{s:appl}) which address this issue.

	This article unfolds as follows. \secref{s:gen} contains our main results---Theorem~\ref{maintheorem} and~\ref{gentheorem}, which provide conditions under which a Lyapunov function of the states is $\Lp 1$-bounded. We establish this $\Lp 1$-boundedness under the assumptions that a certain derived process is a supermartingale outside a compact set, and some more minor conditions.\footnote{It also seems conceivably possible that relaxed Foster-Lyapunov inequalities as in \cite[Condition $\mathbf D(\phi, V, C)$, p. 1356]{ref:doucetal04} arising in the context of subgeometric convergence to a stationary distribution can be employed in the construction of the aforementioned supermartingale; this constitutes future work.} (The supermartingale condition alone is not enough, as pointed out in~\cite{pemantle99}, where the authors establish variants of our results for scalar, possibly non-Markovian processes having increments with bounded $p$-th moments for $p > 2$.) For our results to hold, the underlying process need not be time-homogeneous or Markovian. To wit, in \secref{sec:hybrid} we define a class of hybrid processes that switch between two Markov processes depending on whether they inside or outside a fixed set in the state-space, and demonstrate that although the resulting process may be non-Markovian, our results continue to hold. Connections to optimal stopping problems are drawn in \secref{optstop}, which gives a systematic procedure for verifying our assumptions. In \secref{diffsample} we apply the techniques our techniques to a class of sampled diffusion processes. In addition to the cases considered here, the results in~\secref{s:gen} will be of interest in queueing theory, along the lines of the works~\cite{hastad96, borodin01}. \secref{s:appl} contains some applications of the results in~\secref{s:gen} to stability and robustness of \ifs{}. The classical weak stability questions concerning the existence and uniqueness of invariant measures of \ifs{}, addressed in e.g.,~\cite{diaconisifs, tweedie01, szarek06}, revolve around average contractivity hypotheses of the constituent maps and continuity of the probabilities. In~\secref{s:ifs} we look at stronger stability properties of the \ifs{}, namely, global asymptotic stability almost surely and in expectation, for which we give sufficient conditions. There are no assumptions of global contractivity or memoryless choice of the maps at each iterate; we just require a condition resembling average contractivity in terms of Lyapunov functions with a suitable coupling condition with the Markovian transition probabilities. We mention that although some of the assumptions in~\cite{tweedie01} resemble ours, the conditions needed to establish existence of invariant measures in~\cite{tweedie01} are stronger than what we employ; see~\secref{s:ifs} for a detailed comparison. We also demonstrate in~\secref{s:ifsd} that under mild assumptions, iterated function systems possess strong stability and robustness properties with respect to bounded disturbances. In this subsection the exogenous bounded disturbance is not modeled as a random process.

\subsubsection*{Notations}
	Let $\N \Let  \{1, 2, \ldots\}$, $\Nz \Let  \{0, 1, 2, \ldots\}$, and $\posR \Let  [0, \infty[$. We let $\norm{\cdot}$ denote the standard Euclidean norm on $\R^d$. We let $\bar B_r$ denote the closed Euclidean ball around $0$, i.e., $\bar B_r \Let  \bigl\{y\in\R^d\big| \norm y \le r\bigr\}$. For a vector $v\in\R^d$ let $v\transp$ denote its transpose, and $\norm{v}_{P}$ denote $\sqrt{v\transp P v}$ for a $d\times d$ real matrix $P$. The maximum and minimum of two real numbers $a$ and $b$ is denoted by $a\mx b$ and $a\mn b$, respectively. 

\section{General results}
\label{s:gen}
Before we get into hybrid systems, it will be simpler to follow the arguments if we start by considering a discrete-time Markov chain.

\subsection{Obtaining $\pmb L_1$ Bound using Excursions}
Let $X \Let (X_t)_{t\in\Nz}$ be a discrete time Markov chain with a state space $\sts$. We denote the transition kernel of this chain by $\Pm$, i.e., for every $x\in \sts$, the probability measure $\Pm_x(\cdot)\Let \Pm(x,\cdot)$ determines the law of $X_{t+1}$, conditioned on $X_t=x$. 
At this point we only assume the state space $\sts$ to be any Polish space. 
\begin{assumption}
\label{asssup} 
	There exists a nonnegative function $\varphi:\Nz \times \sts \lra \posR $ satisfying the following.
	\begin{enumerate}[label={\rm (\roman*)}, leftmargin=*, widest=iii, align=right]
		\item There exists a subset $K \subset \sts$ such that the process $(Y_t)_{t\in\Nz}$ defined by $Y_t = \varphi(t,X_t)$ is a supermartingale under $\Pm_{x_0}$, for every $x_0\in \sts \setminus K$  until the first time $X_t$ hits $K$. To wit, if $X_0=x_0 \in \sts\setminus K$ and we define
		\[
			\tauk = \inf\bigl\{ t> 0\;\big|\; X_t \in K\bigr\},
		\]
		then the process $\bigl(Y_{t\wedge \tauk}\bigr)_{t\in\Nz}$ is a supermartingale under $\Pm_{x_0}$. 
% The domain of the function varphi has been modified
		\item There exists a nonnegative measurable real-valued function $V:\sts\lra\rr$ and a positive sequence $(\theta(t))_{t \in \Nz}$ such that
		\[
		\varphi(t, x) \ge V(x) / \theta(t) \quad \text{for all }(t, x)\in\Nz\times\sts,
		\]
		and $C \Let \sum_{t\in\Nz} \theta(t) < \infty$.
		\item $\delta\Let \sup_{x\in K} V(x) < \infty$.		\AssumptionEnd
	\end{enumerate}
\end{assumption}

Our objective is to prove under the above condition (and another minor assumptions) that there exists a bound on $\sup_t\EE_{\xz}\bigl[V(X_t)\bigr]$ depending on $x_0$.

\begin{theorem}
\label{maintheorem}
Consider the setup in Assumption \ref{asssup}, and assume that
\eq\label{whatisbeta}
	\beta\Let  \sup_{x_0\in K} \EE\bigl[\varphi(0,X_1) \indic{\{ X_1 \in \sts\setminus K\}}\; \big|\; X_0=x_0 \bigr] < \infty.
\en
Let $\gamma\Let \sup_{t\in\Nz} \theta(t)$. Then we have
\[
	\sup_{t\in\Nz}\EE_{x_0}\bigl[V(X_t)\bigr] \le C\beta + \delta + \gamma \varphi(0, x_0).
\]
\end{theorem}

In the rest of this section we prove the above theorem. Fix a time $t\in\Nz$, and define two random times
\[
	g_t \Let  \sup\bigl\{s\in\Nz\;\big|\; s\le t,\; X_s\in  K\bigr\} \quad \text{and} \quad h_t \Let  \inf\bigl\{s\in\Nz\;\big|\; s\ge t,\; X_s\in K \bigr\}.
\]
We follow the standard custom of defining supremum over empty sets to be $-\infty$, and the infimum over empty sets to be $+\infty$.

Note that $g_t$ is not a stopping time with respect to the natural filtration generated by the process $X$, although $h_t$ is. The random interval $[g_t,h_t]$ is a singleton if and only if $X_t \in K$. Otherwise, we say that $X_t$ is within an excursion outside $K$. 

Now we have the following decomposition:
\begin{equation}\label{conditioning}
	\EE_{\xz}\bigl[V(X_t)\bigr] =  \EE_{\xz}\bigl[V(X_t)\indic{\{g_t=-\infty\}}\bigr] + \sum_{s=0}^t \EE_{\xz}\bigl[V(X_t)\indic{\{g_t=s\}}\bigr].
\end{equation}
Our first objective is to bound each of the expectations $\EE_{\xz}\bigl[V(X_t)\indic{\{g_t=s\}}\bigr]$.

Before we move on, let us first prove a Lemma which follows readily from Assumption 1. 

\begin{lemma}\label{firstbnd}
	Let $X_0=x_0 \in \sts\setminus K$. Then
	\begin{equation}
	\label{supbnd}
		\EE_{x_0}\bigl[V(X_s)\indic{\{\tau^{\vphantom{T}}_K > s\}}\bigr] \le \varphi(0, x_0)\theta(s) \qquad \text{for } s\in\Nz,
	\end{equation}
	where $(\theta(t))_{t\in\Nz}$ is defined in Assumption \ref{asssup}. 
\end{lemma}
\begin{proof}
	This is a straightforward application of Optional Sampling Theorem (OST) for discrete-time supermartingales. Applying OST for the bounded stopping time $s\wedge \tau^{\vphantom{T}}_K$ to the supermartingale $\bigl(\varphi(t,X_t)\bigr)_{t\in\Nz}$, in view of $\varphi\ge 0$, we have
	\[
	\begin{split}
		\varphi(0, x_0) & \ge \EE_{\xz}\!\bigl[\varphi(s\wedge \tauk, X_{s\wedge \tauk})\bigr] \ge \EE_{\xz}\bigl[\varphi(s,X_s) \indic{\{ \tau^{\vphantom{T}}_K > s \}}\bigr].
	\end{split}
	\]
	Now, by condition (i) in Assumption \ref{asssup}, we can write $\varphi(s,x) \ge  V(x)/\theta(s)$.  Thus, substituting back, one has
	\[
		\varphi(0, x_0) \ge \EE_{\xz}\bigl[V(X_s) \indic{\{ \tau^{\vphantom{T}}_K > s \}}\bigr]/\theta(s).
	\]
	Since $(\theta(t))_{t\in\Nz}$ is positive, we arrive at~\eqref{supbnd}.
\end{proof}

We are ready for the proof of Theorem~\ref{maintheorem}.

\begin{proof}[Proof of Theorem~\ref{maintheorem}]
Let us consider three separate cases:
\medskip

\noindent \textsf{Case 1.} ($-\infty < g_t < t$). In this case $g_t$ can take values $\{0,1,2,\ldots,t-1 \}$. Now, if $s \in \{0, 1, 2, \ldots, t-1\}$, then
\[
	\begin{split}
	\EE_{x_0}\bigl[V(X_t) & \indic{\{g_t=s\}}\bigr ] = \EE_{x_0}\bigl[ V(X_t)\indic{\{X_s\in K\}}\indic{\{X_i\notin K, \; i=s+1,\ldots,t\}}\bigr]\\
	& = \int_{K}\Pm^s(\xz, \mrm dx)\int_{\sts\setminus K} \Pm(x, \mrm dy)\; \EE_y\bigl[V(X_{t-s-1})\indic{\{\tauk > t-s-1\}}\bigr],
	\end{split}
\]
and by~\eqref{firstbnd} it follows that the right-hand side is at most 
\[
	\int_{K} \Pm^s(\xz, \mrm dx)\int_{\sts\setminus K} \Pm(x, \mrm dy)\; \varphi(0, y)\theta(t-s-1).
\]
Thus, one has
\eq\label{eq:inhomo}
	\begin{split}
	\EE_{x_0}\bigl[V(X_t)\indic{\{g_t=s\}}\bigr] & \le \theta(t-s-1) \int_{K} \Pm^s(\xz, \mrm dx)\int_{\sts\setminus K} \Pm(x, \mrm dy)\; \varphi(0, y)\\
	& \le \theta(t-s-1) \sup_{x\in K} \EE_x\bigl[\varphi(0,X_1) \indic{\{X_1 \in \sts\setminus K\}}\bigr]\\
	& = \theta(t-s-1) \beta.
	\end{split}
\en
\noindent \textsf{Case 2.} ($g_t=t$). This is easy, since $X_t\in K$ implies $V(X_t)\le \delta $. Thus
\[
	\EE_{\xz}\bigl[V(X_t) \indic{\{g_t=t\}}\bigr] \le \delta \Pm_{x_0}(X_t \in K) \le \delta.
\]
\noindent{\textsf{Case 3.}} ($g_t=-\infty$). This is the case when the chain started from outside $K$ and has not yet hit $K$, and therefore,
\[
	\EE_{x_0}\bigl[V(X_t)\indic{\{g_t=-\infty\}}\bigr] = \EE_{x_0}\bigl[V(X_t)\indic{\{\tauk > t\}}\bigr] \le \varphi(0,x_0)\theta(t).
\] 
Combining all three cases above, we get the bound:
\begin{equation}
\label{e:key}
	\EE_{x_0}\bigl[V(X_t)\bigr] \le \sum_{s=0}^{t-1} \theta(t-s-1)\beta + \delta + \varphi(0,x_0)\theta(t).
\end{equation}
Maximizing the right-hand side of~\eqref{e:key} over $t$, we arrive at
\[
	\sup_{t\in\Nz} \EE_{x_0}\bigl[V(X_t)\bigr] \le \beta \sum_{s=0}^{\infty} \theta(s) + \delta + \varphi(0,x_0)\sup_{t\in\Nz}\theta(t),
\]
which is the bound stated in the theorem.
\end{proof}

Often it will turn out that $\varphi(t,x)$ is a function $\psi(t, V(x))$ as in the case of the classical Foster-Lyapunov type supermartingales~\cite{meynbk}. In that case $\varphi(t,x)=\epower{\alpha t} V(x)$, for some positive $\alpha$. Thus $\varphi(t,\cdot)$ is a linear function of $V(x)$ for each fixed $t$, with $\theta(t)=\epower{-\alpha t}$, which shows that the sequence $(\theta(t))_{t\in\Nz}$ is summable. See also~\cite{foss04} and the references therein for more general Foster-Lyapunov type conditions. For examples which are not linear see \secref{diffsample}.

\subsection{A Class of Hybrid Processes}
\label{sec:hybrid}
The preceding analysis can be extended for processes which switch their behavior depending on whether the current value is within $K$ or not. They constitute a particularly useful class of controlled processes in which a controller attempts to drive the system into a \emph{target} or \textit{safe set} $K\subset\sts$ whenever the system gets out of $K$ due to its inherent randomness. Below we give a rigorous construction of such a process.

\subsection*{A process $X$ that is $(Y, Z)$-hybrid with respect to $K$}
Consider a pair of Markov chains $(Y,Z)$ where $Y$ is a time-homogeneous Markov chain, and $Z$ is a (possibly) time inhomogeneous Markov chain. We construct a \emph{hybrid discrete-time stochastic process} $X$ by the following recipe: 

Firstly, let the state space for the process be $\mcal{S}^{\mathbb{N}_0}$ along with the natural filtration
\[
\mcal{F}_0 \subseteq \mcal{F}_1 \subseteq \mcal{F}_2 \subseteq \ldots
\]
generated by the coordinate maps.

Secondly, we define the sequence of stopping times $\sigma_0 \Let  \tau_0 \Let  -\infty$ and $\tau_1\le \sigma_1 \le \tau_2 \le \sigma_2 \le \ldots$ by
\[
\begin{split}
	\tau_{i} & \Let  \inf\bigl\{ t > \sigma_{i-1}\;\big|\; X_t \in K \bigr\} \quad\text{and}\\
	\sigma_{i} & \Let  \inf\bigl\{ t > \tau_i\;\big|\; X_t \notin K \bigr\}
\end{split}
\]
for $i \in\Nz$.

Finally, we define the process $X$ as follows: for a measurable $B\subset\sts$,
\begin{alignat*}{2}
	& \text{if } X_t = x, \therex i:\tau_i \le t < \sigma_i,\;&&
	\begin{cases}
		X_t = Y_t,\\
		\Pm\bigl(X_{t+1}\in B \;\big|\,\mcal{F}_t\bigr) = \Pm\bigl(Y_1\in B\;\big|\; Y_0=x\bigr),
	\end{cases}\\
	& \text{if } X_t = x, \therex i:\sigma_i \le t < \tau_{i+1},\;&&
	\begin{cases}
		X_t = Z_t,\\
		\Pm\bigl(X_{t+1}\in B \;\big|\,\mcal{F}_t\bigr) = \Pm\bigl(Z_{t+1-\sigma_i}\in B\;\big|\; Z_{t-\sigma_i}=x\bigr).
	\end{cases}
\end{alignat*}

To wit, the process defined above behaves as the homogeneous chain $Y$ whenever it is inside $K$. Once the process $X$ exits the set $K$, a controller alters the behavior of the chain which, until it enters $K$ again, behaves as a copy of the inhomogeneous chain $Z$ starting from a point outside $K$. The process $X$ is in general non-Markovian due to the possible time inhomogeneity of $Z$. Nevertheless, it is a natural class of examples of switching systems whose Markovian behavior switches in different regions on the state space. We say that $X$ is \emph{$(Y,Z)$-hybrid with respect to $K$}.

The following generalization of Theorem \ref{maintheorem} can be proved along lines of the original proof. The only requirement is a slight modification of the condition \eqref{whatisbeta} which is needed to alter the second inequality in \eqref{eq:inhomo}.  

\begin{theorem}\label{gentheorem}
Consider a stochastic process $X$ that is $(Y,Z)$-hybrid with respect to a measurable $K\subset\sts$ for some homogeneous Markov chain $Y$ and some possibly inhomogeneous Markov chain $Z$.  Suppose Assumption \ref{asssup} holds for the process $Z$ and  
\eq\label{whatisbeta*}
\beta:=\sup_{y_0 \in K}\EE\left[ \varphi(0,Y_1)1_{\{ Y_1 \in \sts\setminus K\}}\; \big|\; Y_0=y_0 \right] < \infty.
\en 
If the process $X$ starts from $x_0\in \sts\setminus K$, we have
\eq\label{genubnd}
	\sup_{t\in\Nz}\EE_{\xz}\bigl[V(X_t)\bigr] \le C\beta + \delta + \gamma \varphi(0, x_0).
\en
\end{theorem}
It is interesting to note that the right side of above bound is a total of individual contributions by the control (for $C$), the choice of $K$ (for $\delta$), and the initial configuration (for $x_0$). We stress that \emph{the conclusion holds even when $X$ is no longer a Markov chain due to the time inhomogeneity of $Z$}. This is important, especially because operator-theoretic bounds like Foster-Lyapunov, or martingale-based bounds do not work in such a case.

\subsection{Connection with Optimal Stopping Problems}
\label{optstop}
Suppose that we are given a Markov chain $Z$ taking values in $\sts$, a function $V:\sts\lra\R$, and a measurable target or safe set $K\subset\sts$. (Alternatively, we may assume that we are given an $\sts$-valued process $X$ that is $(Y,Z)$-hybrid with respect to a measurable $K\subset\sts$.) Our objective is to investigate whether the sequence $\bigl(V(X_t)\bigr)_{t\in\Nz}$ is $\Lp 1$-bounded. To this end one can follow the two-step procedure of first searching for a function $\varphi$ satisfying Assumption \ref{asssup}, followed by an application of Theorem \ref{gentheorem}. A systematic procedure of doing this is given by the following connection with Optimal Stopping problems. 

Let $(\theta(t))_{t\in \Nz}$ be some positive sequence of numbers such that $\sum_{t\in \Nz} \theta(t)$ is finite. Define the \textit{pay-off} or the \textit{reward} function as
\[
h(t,x)= \begin{cases}
V(x)/ \theta(t) &\quad \text{if}\; x\in \sts\setminus K, \; t\in \Nz,\\
0&\quad \text{if}\; x \in K,\; t\in \Nz.
\end{cases}
\] 

Recall that the Optimal Stopping problem \cite[Chapter 1]{PeskirShiryaev} for the process $Z$ and the reward function $h$ defined above consists of finding a stopping time $\tau^*$ such that
\eq\label{ostop}
\EE_{x}\Bigl[ h\bigl(\tau^*\wedge \tauk , Z_{\tau^* \wedge \tauk }\bigr)\Bigr] = \esssup_{\tau} \EE_x\Bigl[ h\bigl(\tau\wedge \tauk , Z_{\tau\wedge\tauk }\bigr)\Bigr],
\en
where $\tauk$ is the hitting time to the set $K$, and $\esssup$ refers to essential supremum over the set of all possible stopping times (see \cite[Chapter 1, Lemma 1.3]{PeskirShiryaev}).

Define the \textit{value function} as
\begin{equation}
\label{valuefunction}
\varphi(n,x_0) \Let \esssup_{\tau \in \stp_n} \EE\bigl[ h(\tau, V(Z_{\tau})) \,\big|\, Z_n = x_0 \bigr],
\end{equation}
where $\stp_n$ is the set of stopping times
\[
\bigl\{(\tau\mx n)\mn\tauk \,\big|\, \text{$\tau$ an arbitrary stopping time} \bigr\}.
\] 

\begin{theorem}
Suppose that the  value function $\varphi(0,x_0)$ is finite for all $x_0\in \sts$, then 
\begin{enumerate}[label={\rm (\roman*)}, leftmargin=*, widest=ii, align=right]
\item $\varphi(t,x_0)$ is finite for all $t\in \Nz$ and 
\[
\varphi(t, x_0) \ge V(x_0)/ \theta(t) \quad \text{for all }(t, \xz)\in\Nz\times (\sts\setmin K).
\]
\item The process $(Y_t)_{t\in\Nz}$ defined by
\[
Y_t \Let \varphi\bigl(t\wedge \tauk , Z_{t\wedge \tauk }\bigr)
\]
is a supermartingale.
\end{enumerate}
\end{theorem}

\begin{proof}
The proof follows from the general theory of optimal stopping. See, for example, \cite[Chapter 4]{CRS}. The sequence of rewards is given by the process $V(Z_{t\wedge \tauk })/\theta(t)$, $t=0,1,2,\ldots$. Applying \cite[Theorem 4.1, p.~66]{CRS} we get
\[
\varphi(n,x_0) = \bigl(V(x_0)/\theta(n)\bigr)\mx\Bigl(\EE\Bigl[ \varphi\bigl(n+1, Z_{(n+1)\wedge \tauk }\bigr)\,\Big|\, Z_{n\wedge \tauk }=x_0  \Bigr] \Bigr).
\]
By considering the first of the two terms in the maximum on the right-hand side above we obtain (i), and (ii) follows from the second.
\end{proof}

In other words, the value function $\varphi(t,x)$ defined in \eqref{valuefunction} satisfies the conditions of Theorem \ref{gentheorem}.

\begin{theorem}
Consider an $\sts$-valued process $X$ that is $(Y,Z)$-hybrid with respect to a measurable $K\subset\sts$ as in \secref{sec:hybrid}. Suppose that for some nonnegative integrable sequence $(\theta(t))_{t \in \Nz}$ the optimal stopping problem \eqref{ostop} has a finite value function $\varphi(t,x_0)$. If additionally condition \eqref{whatisbeta*} is true, then the bound \eqref{genubnd} holds. 
\end{theorem}

Let us remark that the value function, being the \textit{envelope}, is the smallest supermartingale (hence the sharpest bound) that can satisfy Theorem \ref{gentheorem}. Several methods of solving optimal stopping problems in the Markovian setting are available and we refer the reader to \cite{PeskirShiryaev} for a complete review. 

\begin{remark}
	There is a parallel converse result employing standard Foster-Lyapunov techniques for the verification of $f$-ergodicity and $f$-regularity \cite[Chapter 14]{meynbk} of Markov processes. The analysis is based on the functional inequality $\EE[V(X_1)\mid X_0 = x] - V(x) \le -f(x) + b\indic{C}(x)$ for measurable functions $V:\sts\lra[0, \infty]$ and $f:\sts\lra[1, \infty[$, a scalar $b > 0$, and a Borel subset $C$ of $\sts$; \cite[Theorem 14.2.3]{meynbk} asserts that the minimal solution to this inequality, which exists if $C$ is petite (see \cite{meynbk} for precise details), is a ``value function'' given by $G_C(x, f) \Let \EE\bigl[\sum_{t=0}^{\sigma_C} f(X_t)\,\big|\,X_0 = x]$, where $\sigma_C$ is the first hitting-time to $C$. The proof is also based on the existence of a certain supermartingale, and the Markov property is employed crucially.\RemarkEnd
\end{remark}

\subsection{A Class of Sampled Diffusions} \label{diffsample}
In the setting of the process $X$ being $(Y, Z)$-hybrid with respect to a given set $K$, suppose that the state-space for the Markov chains $Y$ and $Z$ is $\rr^d$ and the safe set $K$ is compact. Observe that the only challenge in applying Theorem  \ref{gentheorem} is to find a suitable function $\varphi$ given the Markov chain $Z$ and the function $V$. In applications, a natural choice for the function $V$ is given by square of the Euclidean norm, i.e., $V(x)=\sum_{i=1}^d x_i^2$. For this choice of $V$, we describe below a natural class of examples of Markov chains for which one can construct a $\varphi$ that satisfies part (i) of Assumption \ref{asssup}.	

Consider a diffusion with a possibly time-inhomogeneous drift function, given by the $d$-dimensional stochastic differential equation
\eq\label{multidiff}
{\mrm{d}} X_t = b(t, X_t) \mrm{d}t + \mrm{d} W_t,
\en
where $W_t=(W_t(1), W_t(2), \ldots, W_t(d))$ is a vector of $d$ independent Brownian motions, and $b:\posR\times \R^d\lra\R^d$ is a measurable function.

We will abuse the notations somewhat and construct a function $\varphi: \posR \times \posR \lra \posR$ such that $\bigl(\varphi(t,V(X_t))\bigr)_{t\in\Nz}$ is a supermartingale outside a compact set $K$ and satisfies $\varphi(t, \xi) \ge \xi / \theta(t)$ for some nonnegative sequence $(\theta(t))_{t\in\Nz}$. We define $Z_i=X_{i\wedge \tauk}$ for $i\in\Nz$; $Z$ is the the diffusion sampled at integer time points before hitting $K$. It is clear that $Z$ is a Markov chain such that $\bigl(\varphi(i,V(Z_i))\bigr)_{i\in\Nz}$ is a supermartingale that satisfies the Assumptions \ref{asssup} as long as $\sum_{t \in \Nz} \theta(t) < \infty$. 
	
To construct such a $\varphi$, let us consider a well known family of one-dimensional diffusion, known as the squared Bessel processes (BESQ). This family is indexed by a single nonnegative parameter $\delta \ge 0$ and is described as the unique strong solution of the SDE
\eq\label{besqsde}
\mrm d Y_t = 2\sqrt{Y_t}\,\mrm d\mathfrak{b}_t + \delta\, \mrm dt, \qquad Y_0=y_0 \ge 0,
\en	
where $\mathfrak b \Let (\mathfrak b_t)_{t\in\Nz}$ is a one-dimensional standard Brownian motion. We have the following Lemma:

\begin{lemma}
Let $F:\rr\lra\posR$ be a nonnegative, increasing, and convex function, and fix any terminal time $S > 0$. Define the function
\eq\label{phibesq}
\varphi(t,y) \Let \EE\bigl[F(Y_{S})\,\big|\, Y_t=y \bigr], \quad t \in[0, S],
\en	
where $Y$ solves the SDE \eqref{besqsde}. Then $\varphi$ satisfies the following properties:
\begin{enumerate}[label={\rm (\roman*)}, leftmargin=*, widest=iii]
	\item $\varphi$ is increasing in $y$,\label{besq:varphiincreasing}
	\item $\varphi$ is convex in $y$, and\label{besq:varphiconvex}
	\item $\varphi$ satisfies the partial differential equation\label{besq:pde}
	\eq\label{eqforvp}
		\begin{cases}
			\dfrac{\partial \varphi}{\partial t} + \delta \varphi' + 2y \varphi'' = 0, \quad y > 0, \; t\in(0, S),\\
			\varphi(S, y) = F(y).
		\end{cases}
	\en
\end{enumerate}	
\end{lemma}

Note that $\varphi'$ and $\varphi''$ in the statement of Lemma \ref{besqsde} refers to the first and second derivatives with respect to the second argument of $\varphi$.

\begin{proof}
The proof proceeds by coupling. Let us first show that $\varphi$ is increasing as claimed in \ref{besq:varphiincreasing}. Fix $S > 0$. Consider any two starting points $0 \le x < y$. Construct on the same sample space two copies of BESQ processes $Y^{(1)}$ and $Y^{(2)}$ such that both of them satisfy \eqref{besqsde} with respect to the same Brownian motion $\mathfrak b$ but $Y^{(1)}_0= x$ and $Y^{(2)}_0=y$. It is possible to do this since the SDE \eqref{besqsde} admits a strong solution (see \cite[Chapter 5, Proposition 2.13]{karatshreve}). Hence, by \cite[Chapter 5, Proposition 2.18]{karatshreve}, it follows that $Y^{(1)}_t \le Y^{(2)}_t$ for all $t\ge 0$. Since $F$ is an increasing function, we get
\[
\varphi(t,x)= \EE_x\Bigl[ F\bigl(Y^{(1)}_{S - t}\bigr) \Bigr] \le  \EE_y\Bigl[F\bigl(Y^{(2)}_{S - t}\bigr)\Bigr] = \varphi(t,y).
\] 
This proves that $\varphi$ is increasing in the second argument.

For convexity of $\varphi$ claimed in \ref{besq:varphiconvex}, we use a different coupling. We follow arguments very similar to the one used in the proof of \cite[Theorem 3.1]{hobson98}. Consider three initial points $0 < z < y < x$. And let $\hat X, \hat Y, \hat Z$ be three independent BESQ processes that start from $x,y$, and $z$ respectively. Define the stopping times
\[
\tau_{x}= \inf \Bigl\{ u\;\Big| \; \hat Y_u = \hat X_u  \Bigr\}, \quad \tau_{z}= \inf \Bigl\{ u\;\Big| \; \hat Y_u = \hat Z_u  \Bigr\}. 
\]
Fix a time $t\in[0, S]$, and let $T = S - t$. Define
\[
\sigma = \tau_x \wedge \tau_z \wedge T.
\]
Now, on the event $\sigma=\tau_x$, it follows from symmetry that
\eq\label{coupeq1}
\begin{split}
\EE\Bigl[ \left(\hat X_T - \hat Z_T\right) F\bigl(\hat Y_T\bigr)\indic{\{\sigma=\tau_x\}}\Bigr] &= \EE\Bigl[\left(\hat Y_T - \hat Z_T\right) F\bigl(\hat X_T\bigr)\indic{\{\sigma=\tau_x\}}\Bigr],\\
\EE\Bigl[ \left(\hat X_T - \hat Y_T\right)F\bigl(\hat Z_T\bigr)\indic{\{\sigma=\tau_x\}}\bigr] &= 0.
\end{split}
\en
Similarly, on the event $\sigma=\tau_z$, we have
\eq\label{coupeq2}
\begin{split}
\EE\Bigl[ \left(\hat X_T - \hat Z_T\right) F\bigl(\hat Y_T\bigr)\indic{\{\sigma=\tau_z\}}\Bigr] &= \EE\Bigl[ \left(\hat X_T - \hat Y_T\right) F\bigl(\hat Z_T\bigr)\indic{\{\sigma=\tau_z\}}\Bigr],\\
\EE\Bigl[ \left(\hat Z_T - \hat Y_T\right)F\bigl(\hat X_T\bigr)\indic{\{\sigma=\tau_z\}}\bigr] &= 0.
\end{split}
\en
And finally, when $\sigma=T$, we must have $\hat Z_T < \hat Y_T < \hat X_T$. We use the convexity property of $F$ to get 
\begin{equation}
\label{coupeq3}
\begin{aligned}
\EE\Bigl[ \left(\hat X_T - \hat Z_T\right) F\bigl(\hat Y_T\bigr)\indic{\{\sigma=T\}}\Bigr] & \le \EE\Bigl[ \left(\hat X_T - \hat Y_T\right) F\bigl(\hat Z_T\bigr)\indic{\{\sigma=T\}}\Bigr]\\
	& \quad + \EE\Bigl[ \left(\hat Y_T -\hat Z_T\right) F\bigl(\hat X_T\bigr)\indic{\{\sigma=T\}}\Bigr]. 
\end{aligned}
\end{equation}
Combining the three cases in \eqref{coupeq1}, \eqref{coupeq2}, and \eqref{coupeq3} we get 
\begin{equation}
\label{coupeq4}
\EE\Bigl[\left(\hat X_T - \hat Z_T\right) F\bigl(\hat Y_T\bigr)\Bigr] \le \EE\Bigl[\left(\hat X_T - \hat Y_T\right) F\bigl(\hat Z_T\bigr)\Bigr] + \EE\Bigl[\left(\hat Y_T - \hat Z_T\right) F\bigl(\hat X_T\bigr)\Bigr].
\end{equation}
We now use the fact that $\hat X, \hat Y$, and $\hat Z$ are independent. Also, it is not difficult to see from the SDE \eqref{besqsde} that $\EE_x\bigl[\hat X_T\bigr] - x = \EE_y\bigl[\hat Y_T\bigr] - y= \EE_z\bigl[\hat Z_T\bigr] - z = \delta t$.
Thus, from \eqref{coupeq4} we infer that
\[
(x-z) \varphi(t,y) \le (x-y)\varphi(t,z) + (y-z)\varphi(t,x), \quad \text{for all}\; 0 < z < y < x.
\]
This proves convexity of $\varphi$ in its second argument.

Finally, to see \ref{besq:pde}, it suffices to observe that the equation \eqref{eqforvp} is the classical generator relation for diffusions, for which we refer to \cite[Chapter 5.4]{karatshreve}. The transition density of BESQ processes are smooth and have an explicit representation that satisfy equation \eqref{eqforvp}. The general case can be obtained by differentiating under the integral with respect to $F$.
\end{proof}

Let us return to the multidimensional diffusion given by \eqref{multidiff}. We consider the process $(\zeta_t)_{t\in\Nz}$, where $\zeta_t \Let \varphi\bigl(t, \norm{X_t}^2\bigr)$, and $\varphi$ is the function in \eqref{phibesq}. Note that, since $F$ is nonnegative, so is $\varphi$. Additionally, since $\varphi$ is convex, we have
\[
\varphi(t,\xi) \ge \varphi(t,0) + \varphi'(t,0+)\xi. 
\]
Hence the sequence $(\theta(t))_{t=0}^S$ is given by
\[
\theta(t) = 1 / \varphi'(t,0+), \quad t = 0, 1, \ldots, S.
\]
We have the following Theorem:

\begin{theorem}
Suppose that there exists a compact set $K\subset\R^d$ such that that the drift function $b = (b_1, b_2, \ldots, b_d)$ in the SDE \eqref{multidiff} satisfies the sector condition
\[
\sum_{i=1}^d x_i b_i(t,x) < 0 \quad \text{for } \; (t, x)\in \posR\times (\sts\setmin K).
\]
Fix any terminal time $T > 0$. Define the process $(\zeta_t)_{t\in\Nz} \Let  \bigl(\varphi\bigl(t, \norm{X_t}^2\bigr)\bigr)_{t\in\Nz}$, , where $\varphi$ is the nonnegative, increasing, convex function defined in \eqref{phibesq} with 
\[
F(y)=\norm{y}^2\quad \text{and} \quad \delta = d.
\]
Then, with the set-up as above, the stopped process $\bigl(\zeta_{t\wedge \tauk \wedge T}\bigr)_{t\ge 0}$  is a (local) supermartingale.
% SP: deleted the following bound since the right side is never less than one.
%Moreover, 
%\eq\label{tailbnd}
%\Pm_{x_0}\bigl(\tauk   > s \bigr) \le \inf_{\varphi} \frac{\varphi(0, \norm{x_0}^2)}{\varphi(s,0)}, 
%\en
%where the infimum above is taken over all choices of $\varphi$ for every $T \ge s$.
\end{theorem}

\begin{proof}
Applying It\^o's rule to $(\zeta_t)_{t\in\posR}$, we get
\eq\label{semgle}
\dup \zeta_t = \mrm dM_t + \left[ \frac{\partial \varphi}{\partial t}  + \mathcal{L}\varphi \right]\mrm dt,
\en
where $M \Let (M_t)_{t\in\posR}$ is in general a local martingale ($M$ is a martingale under additional assumptions of boundedness on the first derivative of $\varphi$), and $\mcal{L}$ is the generator of $X$. We compute
\[
\begin{split}
 \frac{\partial \varphi}{\partial t} + \mcal{L} \varphi &=   \frac{\partial \varphi}{\partial t} + \sum_{i=1}^d b_i \frac{\partial \varphi}{\partial x_i} + \frac{1}{2}\sum_{i=1}^d \frac{\partial^2\varphi}{\partial x^2_i} \\
 &= \frac{\partial \varphi}{\partial t} + 2 \varphi' \sum_{i=1}^d b_i x_i + \frac{1}{2}\left[ 2\mrm d \varphi' + \varphi'' \sum_{i=1}^d 4 x_i^2  \right]\\
 &= \frac{\partial \varphi}{\partial t} + \mrm d \varphi' + 2 \left(\sum_i x_i^2\right) \varphi'' + 2 \varphi' \sum_{i=1}^d b_i x_i = 2 \varphi' \sum_i b_i x_i,
\end{split}
\]
where the final equality holds since $\varphi$ satisfies \eqref{eqforvp} at $y = \sum_i x_i^2$.

We know that $\varphi' > 0$ since $\varphi$ is increasing, and, by our assumption, $\sum_i x_i b_i < 0$ whenever $x\not\in K$. Thus, 
\[
 \frac{\partial \varphi}{\partial t} + \mcal{L} \varphi \le 0 \quad \text{for }\; (t, x) \in[0, T]\times (\sts\setmin K).
\] 
Now the claim follows from the semimartingale decomposition given in \eqref{semgle}.
\end{proof}

%To prove the inequality \eqref{tailbnd}, we apply Optional Sampling Theorem \cite[Chapter 1, Theorem 3.22]{karatshreve} to the supermartingale $\Bigl(\varphi\bigl(t\mn\tauk\mn T, \bigl\|X_{t\mn\tauk\mn T}\bigr\|^2\bigr)\Bigr)_{t\in\Nz}$ to obtain
%\[
%\varphi\bigl(0,\norm{x_0}^2\bigr) \ge \EE_{x_0}\Bigl[\varphi\bigl(s\wedge \tauk , \bigl\|X_{s\wedge \tauk }\bigr\|^2\bigr)\Bigr] \ge \EE_{x_0}\Bigl[\varphi\bigl(s, \norm{X_s}^2\bigr) \indic{\{\tauk  > s \}}\Bigr].
%\]
%Since $\varphi$ is increasing in the second argument, we get
%\[
%\varphi\bigl(0,\norm{x_0}^2\bigr) \ge \varphi(s,0)\Pm_{x_0}\bigl(\tauk  > s\bigr).
%\]
%The bound \eqref{tailbnd} follows by rearranging terms above and taking an infimum over the function $\varphi$.

Note that the supermartingale $(\zeta_t)_{t\in\Nz}$ has been defined only for a bounded temporal horizon. Thus, to show that Theorem \ref{gentheorem} holds, some additional uniformity assumptions would be needed.

\medskip

\section{Application to Discrete-Time Randomly Switched Systems}
\label{s:appl}
	In this section we look at several cases of discrete-time randomly switched systems (or, iterated function systems,) in which Theorem~\ref{maintheorem} of~\secref{s:gen} applies and gives useful uniform $\Lp 1$ bounds of Lyapunov functions. In~\secref{s:ifs} we give sufficient conditions for global asymptotic stability almost surely and in $\Lp 1$ of discrete-time randomly switched systems. Assumptions of global contractivity in its standard form or memoryless choice of the maps at each iterate are absent; we simply require a condition resembling average contractivity in terms of Lyapunov functions with a suitable coupling condition with the Markovian transition probabilities. In~\secref{s:ifsd} we demonstrate that under mild hypotheses iterated function systems possess strong stability and robustness properties with respect to bounded disturbances that are not modelled as random processes.\footnote{Recall the following notation: We let $\ClassK$ denote the collection of strictly increasing continuous functions $\alpha:\posR\lra\posR$ such that $\alpha(0) = 0$; we say that a function $\alpha$ belongs to class-$\ClassKinfty$ if $\alpha\in\ClassK$ and $\lim_{r\to\infty}\alpha(r) = \infty$. A function $\beta:\posR\times\Nz\lra\posR$ belongs to class-$\ClassKL$ if $\beta(\cdot, n)\in\ClassK$ for a fixed $n\in\Nz$, and if $\beta(r, n)\to 0$ as $n\to\infty$ for fixed $r\in\posR$. Recall that a function $f:\R^d\lra\R^d$ is locally Lipschitz continuous if for every $\xz\in\R^d$ and open set $O$ containing $\xz$, there exists a constant $L > 0$ such that $\norm{f(x) - f(\xz)} \le L\norm{x-\xz}$ whenever $x\in O$.}

	\subsection{Stability of Discrete-Time Randomly Switched Systems}
	\label{s:ifs}
		Consider the system 
		\begin{equation}
		\label{e:sys}
			X_{t+1} = f_{\sigma_t}(X_t),\qquad X_0 = \xz, \quad t\in\Nz.
		\end{equation}
		Here $\sigma:\Nz\lra\mcal P \Let  \{1, \ldots, \cardP\}$ is a discrete-time random process, the map $f_i:\R^d\lra\R^d$ is continuous and locally Lipschitz, and there are points $x_i^\star\in\R^d$ such that $f_i(x_i^\star) = 0$ for each $i\in\PSet$. The initial condition of the system $\xz\in\R^d$ is assumed to be known. Our objective is to study stability properties of this system by extracting certain nonnegative supermartingales.

		The system~\eqref{e:sys} can be viewed as an iterated function system: $X_{t+1} = f_{\sigma_t}\circ\cdots\circ f_{\sigma_1}\circ f_{\sigma_0}(\xz)$. Varying the point $\xz$ but keeping the same maps leads to a family of Markov chains initialized from different initial conditions. The article~\cite{diaconisifs} treats basic results on convergence and stationarity properties of such systems with the process $(\sigma_t)_{t\in\Nz}$ being a sequence of independent and identically distributed random variables taking values in $\PSet$, and each map $f_i$ is a contraction. These results were generalized in~\cite{tweedie01} with the aid of Foster-Lyapunov arguments.

		The analysis carried out in~\cite{tweedie01} requires a Polish state-space, and employs the following three principal assumptions: (a) the maps are non-separating on an average, i.e., the average separation of the Markov chains initialized at different points is nondecreasing over time; (b) there exists a set $C$ such that the Markov chains started at different initial conditions contract after the set $C$ is reached; and (c) there exists a measurable real-valued function $V\ge 1$, bounded on $C$, and satisfying a Foster-Lyapunov drift condition $QV(x) \le \lambda V(x) + b\indic{C}(x)$ for some $\lambda\in\:]0, 1[$ and $b < \infty$, where $Q$ is the transition kernel. Under these conditions the authors establish the existence and uniqueness of an invariant measure which is also globally attractive, and the convergence to this measure is exponential. In particular, this showed that the main results of~\cite{diaconisifs}, which are primarily related to existence and uniqueness of invariant probability measures, continue to hold if the contractivity hypotheses on the family $\{f_i\}_{i\in\PSet}$ are relaxed. In this subsection we look at stronger properties, namely, $\Lp 1$ boundedness and stability, and almost sure stability of the system~\eqref{e:sys} under Assumption~\ref{asssup}. No contractivity inside a compact set is needed to establish existence of an invariant measure under Assumption~\ref{asssup}.

		\begin{assumption}
		\label{a:sigma}
			The process $(\sigma_t)_{t\in\Nz}$ is an irreducible Markov chain with initial probability distribution $\pi^\circ$ and a transition matrix $P \Let  [p_{ij}]_{\cardP\times\cardP}$.\AssumptionEnd
		\end{assumption}
		It is immediately clear that the discrete-time process $(\sigma_t, X_t)_{t\in\Nz}$, taking values in the Borel space $\PSet\times\R^d$, is Markovian under Assumption~\ref{a:sigma}. The corresponding transition kernel is given by
		\begin{align*}
			Q\bigl((i, x), \mcal P'\times B\bigr) = \textstyle{\sum_{j\in\mcal P'} p_{ij}\indic{B}\bigl(f_j(x)\bigr)}\quad & \text{for }\mcal P'\subset\mcal P,B\text{ a Borel subset of } \R^d,\\
			& \text{and }(i, x)\in\PSet\times\R^d.
		\end{align*}

		Our basic analysis tool is a family of Lyapunov functions, one for each subsystem, and at different times we shall impose the following two distinct sets of hypotheses on them.\footnote{It will be useful to recall here that the deterministic system $x_{t+1} = f_i(x_{t}), \;t\in\Nz,$ with initial condition $\xz$ is said to be \emph{globally asymptotically stable} (in the sense of Lyapunov) if (a) for every $\eps > 0$ there exists a $\delta > 0$ such that $\norm{\xz-x_i^\star} < \delta$ implies $\norm{x_t-x_i^\star} < \eps$ for all $t\in\Nz$, and (b) for every $r, \eps' > 0$ there exists a $T > 0$ such that $\norm{\xz-x_i^\star} < r$ implies $\norm{x_t-x_i^\star} < \eps$ for all $t > T$. The condition (a) goes by the name of Lyapunov stability of the dynamical system (or of the corresponding equilibrium point $x_i^\star$), and (b) is the standard notion of global asymptotic convergence to $x_i^\star$.} 

		\begin{assumption}
		\label{a:V}
			There exist a family $\{V_i\}_{i\in\mcal P}$ of nonnegative measurable functions on $\R^d$, functions $\alpha_1, \alpha_2\in\ClassK$, numbers $\lambda_\circ \in\;]0, 1[$, $r > 0$ and $\mu > 1$, such that
			\begin{enumerate}[label=(V\arabic*), align=right, leftmargin=*]
				\item $\alpha_1(\norm{x-x_i^\star}) \le V_i(x) \le \alpha_2(\norm{x-x_i^\star})\quad$ for all $x$ and $i$,
				\item $V_i(x) \le \mu V_j(x)\quad$ whenever $\norm{x} > r$, for all $i, j$, and
				\item $V_i(f_i(x)) \le \lambda_\circ V_i(x)\quad$ for all $x$ and $i$.\AssumptionEnd
			\end{enumerate}
		\end{assumption}

		\begin{assumption}
		\label{a:Vprime}
			There exist a family $\{V_i\}_{i\in\mcal P}$ of nonnegative measurable functions on $\R^d$, functions $\alpha_1, \alpha_2\in\ClassK$, a matrix $[\lambda_{ij}]_{\cardP\times\cardP}$ with nonnegative entries, and numbers $r > 0$, $\mu > 1$, such that (V1)-(V2) of Assumption~\ref{a:V} hold, and
			\begin{enumerate}[label=(V3$'$), align=left, leftmargin=*]
				\item $V_i(f_j(x)) \le \lambda_{ij} V_i(x)\quad$ for all $x$ and $i, j$.\AssumptionEnd
			\end{enumerate}
		\end{assumption}

		The condition (V1) in Assumption~\ref{a:V} is standard in deterministic system theory literature, ensuring, in particular, positive definiteness of each $V_i$. (V2) stipulates that outside $\bar B_r$ the functions $\{V_i\}_{i\in\PSet}$ are linearly comparable to each other. %This is possible if, for instance, $V_1(x) = \norm{x-x_1^\star}^4_{P_1}$ and $V_2(x) = \norm{x-x_2^\star}^2_{P_2} + \norm{x-x_2^\star}^4_{P_3}$ for some $x_1^\star, x_2^\star\in\R^d$ and positive definite matrices $P_i,\; i=1, 2, 3$. 
		The conditions (V1) and (V3) together imply that each subsystem is globally asymptotically stable, with sufficient stability margin---the smaller the number $\lambda_\circ$, the greater is the stability margin. In fact, standard converse Lyapunov theorems show that (V1) and (V3) are necessary and sufficient conditions for each subsystem to be globally asymptotically stable. The only difference between Assumptions \ref{a:V} and \ref{a:Vprime} is that the latter keeps track of how each Lyapunov function evolves along trajectories of every subsystem.%Note that $\lambda_{ii}$ is not required to be less than $1$ for each $i$, which corresponds to allowing unstable subsystems.

		Let us define %$\hat p$ and $\tilde p$ to be the maximum diagonal and off-diagonal entry of the matrix $P$, respectively, i.e., 
		$\displaystyle{\hat p \Let  \max_{i\in\cardP} p_{ii}}$ and $\displaystyle{\tilde p \Let  \max_{i, j\in\PSet, i\neq j} p_{ij}}$.

		\begin{proposition}
		\label{p:supmartss}
			Consider the system~\eqref{e:sys}, and suppose that either of the following two conditions holds:
			\begin{enumerate}[label=\emph{(S\arabic*)}, align=right, leftmargin=*]
				\item Assumptions~\ref{a:sigma} and~\ref{a:V} hold, and $\lambda_\circ(\hat p + \mu\tilde p) < 1$.
				\item Assumptions~\ref{a:sigma} and~\ref{a:Vprime} hold, and $\textstyle{\mu\cdot\left(\max_{i\in\PSet} \sum_{j\in\PSet} p_{ij}\lambda_{ji}\right) < 1}$.
			\end{enumerate}
			Let $\tau_r \Let  \inf\bigl\{t\in\Nz\big|\norm{X_t} \le r\bigr\}$ and $V_i'(x) \Let  V_i(x)\indic{\R^d\setminus \bar B_r}(x)$. Suppose that $\norm{\xz} > r$. Then there exists $\alpha > 0$ such that the process $\bigl(\epower{\alpha (t\mn\tau_r)} V'_{\sigma_{t\mn\tau_r}}(X_{t\mn\tau_r})\bigr)_{t\in\Nz}$ is a nonnegative supermartingale.% In particular, there exists a constant $C > 0$ such that $\sup_{t\in\Nz} \Ex{V_{\sigma_t}(X_t)} < C$.
		\end{proposition}

		%Note that the condition (S1) of Proposition~\ref{p:supmartss} exploits less of the Markovian structure of $(\sigma_t)_{t\in\Nz}$ than the condition (S2)---just two elements of the transition matrix $[p_{ij}]_{\cardP\times\cardP}$ are involved in (S1), whereas all the elements of this matrix are involved in (S2). On the one hand, the inequality in (S1) is clearly conservative, since we employ the uniform bounds of the transition probability matrix and stability margins of the subsystems. None of the subsystems is permitted to be unstable or expansive, for $\lambda_\circ$ is assumed to be strictly less than $1$. On the other hand, the inequality in (S2) is a weighted sum, and unstable or expansive subsystems are permissible as long as the corresponding weights are small. Finally, it is necessary to keep track of how each Lyapunov function evolves along the trajectories of every other subsystem if unstable or expansive subsystems are allowed, which is why we need all the entries of the matrix $[\lambda_{ij}]_{\cardP\times\cardP}$ in (S2), unlike (S1) where we employ uniform stability bounds captured by $\lambda_\circ$.

		\begin{corollary}
		\label{c:boundss}
			Consider the system~\eqref{e:sys}, and assume that the hypotheses of Proposition~\ref{p:supmartss} hold. Then there exists a constant $c > 0$ such that $\displaystyle{\sup_{t\in\Nz} \Ex{\alpha_1(\norm{X_t})} < c}$.
		\end{corollary}

		It is possible to derive simple conditions for stability of the system~\eqref{e:sys} from Proposition~\ref{p:supmartss}. To this end we briefly recall two standard stability concepts.
		\begin{defn}
		\label{d:gasas}
			If $\ker(f_i - \id) = \{0\}$ for each $i\in\PSet$, the system~\eqref{e:sys} is said to be 
			\begin{itemize}[label=$\circ$, leftmargin=*]
				\item \emph{globally asymptotically stable almost surely} if
				\begin{enumerate}[label=(AS\arabic*), align=right, leftmargin=*]
					\item $\displaystyle{\msf P\Bigl(\fa \eps > 0\;\; \therex \delta > 0\text{ s.t.\ } \sup_{t\in\Nz} \norm{X_t} < \eps\text{ whenever } \norm{\xz} < \delta\Bigr) = 1}$,
					\item $\displaystyle{\msf P\Bigl(\fa r, \eps' > 0\;\; \therex T > 0\text{ s.t.\ } \sup_{\Nz\ni t > T} \norm{X_t} < \eps'\text{ whenever } \norm{\xz} < r\Bigr) = 1}$;
				\end{enumerate}
				\item \emph{$\alpha$-stable in $\Lp 1$} for some $\alpha \in\ClassK$ if
				\begin{enumerate}[label=(SM\arabic*), align=right, leftmargin=*]
					\item $\displaystyle{\fa \eps > 0\;\; \therex \delta > 0\text{ s.t.\ } \sup_{t\in\Nz} \Ex{\alpha(\norm{X_t})} < \eps\text{ whenever } \norm{\xz} < \delta}$,
					\item $\displaystyle{\fa r, \eps' > 0\;\; \therex T > 0\text{ s.t.\ } \sup_{\Nz\ni t > T} \Ex{\alpha(\norm{X_t})} < \eps'\text{ whenever } \norm{\xz} < r}$.\DefEnd
				\end{enumerate}
			\end{itemize}
		\end{defn}

		\begin{corollary}
		\label{c:gasas}
			Suppose that $\ker(f_i - \id) = \{0\}$ for each $i\in\PSet$, and that either of the hypotheses \emph{(S1)} and \emph{(S2)} of Proposition~\ref{p:supmartss} holds with $r = 0$. Then
			\begin{itemize}[label=$\circ$, leftmargin=*]
				\item there exists $\alpha > 0$ such that $\lim_{t\to\infty}\msf E\bigl[\epower{\alpha t}V_{\sigma_{t}}(X_{t})\bigr] = 0$, and
				\item the system~\eqref{e:sys} is globally asymptotically stable almost surely and $\alpha_1$-stable in $\Lp 1$ in the sense of Definition~\ref{d:gasas}.
			\end{itemize}
		\end{corollary}

		The proofs of Proposition~\ref{p:supmartss}, Corollary~\ref{c:boundss} and Corollary~\ref{c:gasas} are given after the following simple Lemma; the crude estimate asserted in it resembles the distribution of a Binomial random variable, except that we have $\hat p + \tilde p \ge 1$. For $t\in\N$ let the random variable $N_t$ denote the number of times the state of the Markov chain changes on the period of length $t$ starting from $0$, i.e., $N_t \Let  \sum_{i=1}^{t} \indic{\{\sigma_{i-1} \neq \sigma_{i}\}}$.

		\begin{lemma}
		\label{l:sigmalaw}
			Under Assumption~\ref{a:sigma} we have for $s < t$, $s, t\in\Nz$,
			\[
				\msf P\bigl(N_t - N_s = k\big|\sigma_s\bigr) \le 
				\begin{cases}
					\displaystyle{\left(\binom{t-s}{k}\hat p^{(t-s-k)}\tilde p^k\right)\mn 1} \quad & \text{if $k=0, 1,\ldots, t-s$},\\
					0 & \text{else}.
				\end{cases}
			\]
		\end{lemma}
		\begin{proof}
			Fix $s < t$, $s, t\in\Nz$, and let $\eta_k(s, t) \Let  \msf P\bigl(N_t - N_s = k\big|\sigma_s\bigr)$. Then by the Markov property, for $k = 0, 1, \ldots, t-s$,
			\begin{align*}
				\eta_k(s, t) %& = \msf P\bigl(N_t - N_s = k\big|N_{t-1} - N_s = k, \sigma_s\bigr)\msf P\bigl(N_{t-1} - N_s = k\big|\sigma_s\bigr)\\
				%& \qquad + \msf P\bigl(N_t - N_s =k\big|N_{t-1} - N_s = k-1, \sigma_s\bigr)\msf P\bigl(N_{t-1} - N_s = k-1\big|\sigma_s\bigr)\\
				& = \eta_k(s, t-1) \msf P\bigl(N_t - N_s = k\big|N_{t-1} - N_s = k, \sigma_s\bigr)\\
				& \qquad + \eta_{k-1}(s, t-1)\msf P\bigl(N_t - N_s = k\big|N_{t-1} - N_s = k-1, \sigma_s\bigr)\\
				& \le \hat p\eta_k(s, t-1) + \tilde p\eta_{k-1}(s, t-1).
			\end{align*}
			The set of initial conditions $\eta_i(s, t) = 0$ for all $i\ge t-s$, follow from the trivial observation that there cannot be more than $t-s$ changes of $\sigma$ on a period of length $t-s$. This gives a well-defined set of recursive equations, and a standard induction argument shows that $\eta_k(s, t) \le \binom{t-s}{k}\hat p^{(t-s-k)}\tilde p^{k}$. This proves the assertion.
		\end{proof}

		%The smaller the difference between the maximal off-diagonal and diagonal elements of the transition matrix $[p_{ij}]_{\cardP\times\cardP}$, the tighter is the bound in Lemma~\ref{l:sigmalaw}. The bound also becomes tighter as the difference $(t-s)$ increases.

		\begin{proof}[Proof of Proposition~\ref{p:supmartss}]
			First we look at the assertion under the condition (S1). Fix $s < t$, $s, t\in\Nz$. Given $(\sigma_{s\mn\tau_r}, X_{s\mn\tau_r})$, from (V3) we get $V'_{\sigma_{s\mn\tau_r}}\bigl(X_{(s+1)\mn\tau_r}\bigr) \le \lambda_\circ V'_{\sigma_{s\mn\tau_r}}\bigl(X_{s\mn\tau_r}\bigr)$, and if $\sigma_{s+1}\neq \sigma_s$, we employ (V2) to get $V'_{\sigma_{(s+1)\mn\tau_r}}\bigl(X_{(s+1)\mn\tau_r}\bigr) \le \mu V'_{\sigma_{s\mn\tau_r}}\bigl(X_{(s+1)\mn\tau_r}\bigr)$. Therefore,
			\begin{alignat*}{2}
				& V'_{\sigma_{(s+1)\mn\tau_r}}\bigl(X_{(s+1)\mn\tau_r}\bigr) \le \mu\lambda_\circ V'_{\sigma_{s\mn\tau_r}}\bigl(X_{s\mn\tau_r}\bigr) &\quad& \text{if }\sigma_{(s+1)\mn\tau_r}\neq\sigma_{s\mn\tau_r},\quad\text{and}\\
				& V'_{\sigma_{(s+1)\mn\tau_r}}\bigl(X_{(s+1)\mn\tau_r}\bigr) \le \lambda_\circ V'_{\sigma_{s\mn\tau_r}}\bigl(X_{s\mn\tau_r}\bigr) && \text{otherwise}.
			\end{alignat*}
			Iterating this procedure we arrive at the pathwise inequality
			\begin{equation}
			\label{e:pathwise}
				V'_{\sigma_{t\mn\tau_r}}\bigl(X_{t\mn\tau_r}\bigr) \le \mu^{N_{t\mn\tau_r} - N_{s\mn\tau_r}} \lambda_\circ^{t\mn\tau_r-s\mn\tau_r} V'_{\sigma_{s\mn\tau_r}}\bigl(X_{s\mn\tau_r}\bigr).
			\end{equation}
			Since $s\mn\tau_r = t\mn s \mn\tau_r$, and $t\mn\tau_r$ is measurable with respect to $\sigalg_{t\mn s\mn\tau_r}$, we invoke the Markov property of $(\sigma_t, X_t)_{t\in\Nz}$to arrive at
			\[
			\begin{aligned}
				\msf E\Bigl[V'_{\sigma_{t\mn\tau_r}}\bigl(X_{t\mn\tau_r}\bigr) & \Big|(\sigma_{s\mn\tau_r}, X_{s\mn\tau_r})\Bigr]\\
					& \le V'_{\sigma_{s\mn\tau_r}}\bigl(X_{s\mn\tau_r}\bigr) \lambda_\circ^{t\mn\tau_r-s\mn\tau_r} \msf E\Bigl[\mu^{N_{t\mn\tau_r} - N_{s\mn\tau_r}}\Big|(\sigma_{s\mn\tau_r}, X_{s\mn\tau_r})\Bigr].
			\end{aligned}
			\]
			We now apply the estimate in Lemma~\ref{l:sigmalaw} to get $\msf E\Bigl[\mu^{N_{t\mn\tau_r} - N_{s\mn\tau_r}}\Big|(\sigma_{s\mn\tau_r}, X_{s\mn\tau_r})\Bigr] \le \sum_{k=0}^{t\mn\tau_r-s\mn\tau_r} \binom{t\mn\tau_r-s\mn\tau_r}{k} \hat p^{(t\mn\tau_r-s\mn\tau_r-k)} \tilde p^k \mu^{k} = \bigl(\hat p + \mu\tilde p\bigr)^{t\mn\tau_r-s\mn\tau_r}$,
			and this leads to
			\[
				\msf E\Bigl[V'_{\sigma_{t\mn\tau_r}}\bigl(X_{t\mn\tau_r}\bigr)\Big|(\sigma_{s\mn\tau_r}, X_{s\mn\tau_r})\Bigr] \le V'_{\sigma_{s\mn\tau_r}}\bigl(X_{s\mn\tau_r}\bigr) \bigl(\lambda_\circ (\hat p + \mu\tilde p)\bigr)^{t\mn\tau_r-s\mn\tau_r}.
			\]
			Since $\lambda_\circ(\hat p + \mu\tilde p) < 1$, letting $\alpha' \Let  \lambda_\circ(\hat p + \mu\tilde p)\epower{\alpha} < 1$, the above inequality gives
			\begin{equation}
			\label{e:supmartineq}
			\begin{aligned}
				\msf E\Bigl[\epower{\alpha (t\mn\tau_r-s\mn\tau_r)}V'_{\sigma_{t\mn\tau_r}} & \bigl(X_{t\mn\tau_r}\bigr)\Big|(\sigma_{s\mn\tau_r}, X_{s\mn\tau_r})\Bigr]\\
					& \le V'_{\sigma_{s\mn\tau_r}}\bigl(X_{s\mn\tau_r}\bigr) (\alpha')^{t\mn\tau_r-s\mn\tau_r} \le V'_{\sigma_{s\mn\tau_r}}\bigl(X_{s\mn\tau_r}\bigr).
			\end{aligned}
			\end{equation}
			This shows that $\bigl(\epower{\alpha (t\mn\tau_r)}V'_{\sigma_{t\mn\tau_r}}\bigl(X_{t\mn\tau_r}\bigr)\bigr)_{t\in\Nz}$ is a nonnegative supermartingale.

			Let us now look at the assertion of the Proposition under the condition (S2). Fix $t\in\Nz$. Then from (V3$'$), $V'_j\bigl(f_{\sigma_{t\mn\tau_r}}\bigl(X_{t\mn\tau_r}\bigr)\bigr) \le \lambda_{j\sigma_{t\mn\tau_r}}V'_j\bigl(X_{t\mn\tau_r}\bigr)$ for all $j\in\PSet$, and by (V2), $V'_{\sigma_{(t+1)\mn\tau_r}}\bigl(f_{\sigma_{t\mn\tau_r}}\bigl(X_{t\mn\tau_r}\bigr)\bigr) \le \lambda_{\sigma_{(t+1)\mn\tau_r}\sigma_{t\mn\tau_r}} V'_{\sigma_{(t+1)\mn\tau_r}}\bigl(X_{t\mn\tau_r}\bigr) \le \mu\lambda_{\sigma_{(t+1)\mn\tau_r}\sigma_{t\mn\tau_r}} V'_{\sigma_{t\mn\tau_r}}\bigl(X_{t\mn\tau_r}\bigr)$.  This leads to
			\[%begin{equation}
			\begin{aligned}
			%\label{e:supmartineq2}
				\msf E\Bigl[V'_{\sigma_{(t+1)\mn\tau_r}}\bigl(X_{(t+1)\mn\tau_r}\bigr) \Big|(\sigma_{t\mn\tau_r}, X_{t\mn\tau_r})\Bigr] %\le \mu\Biggl(\sum_{j\in\PSet} p_{\sigma_{t\mn\tau_r}j}\lambda_{j\sigma_{t\mn\tau_r}}\Biggr) V'_{\sigma_{t\mn\tau_r}}(X_{t\mn\tau_r})\\
				& \le \mu\Biggl(\max_{i\in\PSet}\sum_{j\in\PSet} p_{ij}\lambda_{ji}\Biggr) V'_{\sigma_{t\mn\tau_r}}\bigl(X_{t\mn\tau_r}\bigr).
			\end{aligned}
			\]%end{equation}
			Since by hypothesis there exists $\alpha > 0$ such that $\mu\left(\max_{i\in\PSet}\sum_{j\in\PSet} p_{ij}\lambda_{ji}\right)\epower{\alpha} < 1$, the last inequality shows immediately that $\bigl(\epower{\alpha (t\mn\tau_r)}V'_{\sigma_{t\mn\tau_r}}\bigl(X_{t\mn\tau_r}\bigr)\bigr)_{t\in\Nz}$ is a supermartingale. This concludes the proof.
		\end{proof}

		%\begin{remark}
		%	It is possible to modify the above proof of Proposition~\ref{p:supmartss} to directly verify that the process $\bigl(\epower{\alpha (t\mn\tau_r)}V'_{\sigma_{t\mn\tau_r}}\bigl(X_{t\mn\tau_r}\bigr)\bigr)_{t\in\Nz}$ is a uniformly integrable supermartingale. Indeed, minor modifications in the proof are needed to demonstrate that the process $\Bigl(\bigl(\epower{\alpha (t\mn\tau_r)}V'_{\sigma_{t\mn\tau_r}}\bigl(X_{t\mn\tau_r}\bigr)\bigr)^{1+\delta}\Big)_{t\in\Nz}$ is a supermartingale for some $\delta > 0$ small enough; the key observation is that the maps $]0, \infty[\;\ni y \mapsto \lambda_\circ^{1+y}\bigl(\hat p + \tilde p\mu^{1+y}\bigr)\in\;]0, \infty[$ and $]0, \infty[\;\ni y\mapsto \mu^{1+y}\Bigl(\sum_{j\in\PSet} p_{ij}\lambda_{ji}^{1+y}\Bigr)\in\;]0, \infty[$ are continuous.\RemarkEnd
		%\end{remark}

		\begin{proof}[Proof of Corollary~\ref{c:boundss}]
			First observe that since each map $f_i$ is locally Lipschitz, the diameter of the set $D_i \Let  \bigl\{f_i(x)\big| x\in\bar B_r\bigr\}$ is finite, and since $\PSet$ is finite, so is the diameter of $\bigcup_{i\in\PSet} D_i$. Therefore, if $Q$ is the transition kernel of the Markov process $(\sigma_t, X_t)_{t\in\Nz}$, then employing (V1) and the fact that $f_i$ is locally Lipschitz for each $i$, we arrive at%for $\norm{\xz} < r$,
			\begin{align*}
				\msf E\Bigl[V_{\sigma_1}(X_1)& \indic{\{X_1\in\R^d\setminus\bar B_r\}}\Big| (\sigma_0, X_0) = (i, \xz)\Bigr] %= \int_{\PSet\times(\R^d\setminus \bar B_r)} Q\bigl((i, \xz), j\otimes \mrm dy\bigr) V_j(y)\\
				 = \sum_{j\in\PSet}p_{ij}\indic{\R^d\setminus\bar B_r}(f_j(\xz))V_j(f_j(\xz))\\
				 & \le \sum_{j\in\PSet} p_{ij}\indic{\R^d\setminus\bar B_r}(f_j(\xz))\alpha_2(\norm{f_j(\xz)}) \le \sum_{j\in\PSet} p_{ij} L\norm{\xz} < Lr < \infty
				%& < \infty.
			\end{align*}
			for $\norm{\xz} < r$, where $L$ is such that $\sup_{j\in\PSet, y\in\bar B_r}\norm{f_j(y)} \le L\norm{y}$. This shows that condition \textrm{\ref{whatisbeta}} of Theorem~\ref{maintheorem} holds under our hypotheses, and by Proposition~\ref{p:supmartss} we know that there exists $\alpha > 0$ such that $\Bigl(\epower{\alpha(t\mn\tau_r)}V_{\sigma_{t\mn\tau_r}}\bigl(X_{t\mn\tau_r}\bigr)\indic{\R^d\setminus \bar B_r}\bigl(X_{t\mn\tau_r}\bigr)\Bigr)_{t\in\Nz}$ is a supermartingale. Theorem~\ref{maintheorem} now guarantees the existence of a constant $C' > 0$ such that $\sup_{t\in\Nz} \Ex{V_{\sigma_t}(X_t)\indic{\R^d\setminus \bar B_r}(X_t)} \le C'$, and finally, from (V1) it follows that there exists a constant $c > 0$ such that $\sup_{t\in\Nz} \Ex{\alpha_1(\norm{X_t})} \le c < \infty$, as asserted.
		\end{proof}

		\begin{proof}[Proof of Corollary~\ref{c:gasas}]
			We prove almost sure global asymptotic stability and $\alpha_1$ stability in $\Lp 1$ of~\eqref{e:sys} under the condition (S1) of Proposition~\ref{p:supmartss}; the proofs under (S2) are similar.

			First observe that since $\ker(f_i - \id) = \{0\}$ for each $i\in\PSet$, i.e., $0$ is the equilibrium point of each individual subsystem, $\msf P_{\xz}\bigl(\tau_{\{0\}} < \infty\bigr) = 0$ for $\xz\neq 0$, where $\tau_{\{0\}}$ is the first time that the process $(X_t)_{t\in\Nz}$ hits $\{0\}$. Indeed, since $\ker(f_i - \id) = \{0\}$ for each $i\in\PSet$ and $\xz\neq 0$ we have $Q\bigl((i, \xz), \mcal P\times \{0\}\bigr) = \sum_{j\in\PSet} p_{ij}\indic{\{0\}}(f_j(\xz)) = 0$, which shows that $Q^n\bigl((i, \xz), \mcal P\times\{0\}\bigr) = 0$ whenever $\xz\neq 0$. The observation now follows from $\msf P_{\xz}\bigl(\tau_{\{0\}} < \infty\bigr) = \msf P_{\xz}\bigl(\bigcup_{n\in\N}\bigl\{\tau_{\{0\}} = n\bigr\}\bigr) \le \sum_{n\in\N}\msf P_{\xz}\bigl(\tau_{\{0\}} = n\bigr)$. Therefore, with $\tau_{\{0\}} = \tau_r = \infty$, proceeding as in the proof of Proposition~\ref{p:supmartss} above, one can show that $\bigl(\epower{\alpha t}V_{\sigma_{t}}\bigl(X_{t}\bigr)\bigr)_{t\in\Nz}$ is a supermartingale for some $\alpha > 0$. In particular, With $s = 0$ and $\tau_r = \infty$ in~\eqref{e:supmartineq}, we apply (V1) to arrive at $\lim_{t\to\infty} \msf E\bigl[\epower{\alpha t} V_{\sigma_{t}}\bigl(X_{t}\bigr)\bigr] = \lim_{t\to\infty}\msf E\Bigl[\msf E\bigl[\epower{\alpha t} V_{\sigma_{t}}\bigl(X_{t}\bigr)\big|(\sigma_0, x_0)\bigr]\Bigr] \le \lim_{t\to\infty}\alpha_2(\norm{\xz})(\alpha')^t = 0$. Standard supermartingale convergence results and the definition of $\tau_{\{0\}}$ imply that $\msf P\Bigl(\lim_{t\to\infty} V_{\sigma_t}(X_t) = 0\Bigr) = 1$. With $s = 0$ and $\tau_r = \tau_{\{0\}} = \infty$, the pathwise inequality~\eqref{e:pathwise} in conjunction with (V1) give $V_{\sigma_t}(X_t) \le \alpha_2(\norm{\xz}) \mu^{N_t} \lambda_\circ^t$. The foregoing inequality implies that for almost every sample path $(\sigma_t, X_t')_{t\in\Nz}$ corresponding to initial condition $X_0 = \xz'$ with $\norm{\xz'} < \norm{\xz}$, one has 
			\[
				\lim_{t\to\infty} V_{\sigma_t}(X_t') \le \lim_{t\to\infty}\alpha_2(\norm{\xz'})\mu^{N_t} \lambda_\circ^t \le \lim_{t\to\infty}\alpha_2(\norm{\xz})\mu^{N_t}\lambda_\circ^t = 0,
			\]
			which proves (AS2). Since the family $\{f_i\}_{i\in\PSet}$ is finite, and each $f_i$ is locally Lipschitz, there exists $L > 0$ such that $\sup_{i\in\PSet} \norm{f_i(x)} \le L\norm{x}$ whenever $\norm{x} \le 1$. Fix $\eps > 0$. By (AS2) we know that for almost all sample paths there exists a constant $T > 0$ such that $\sup_{t\ge T}\norm{X_t} < \eps$ whenever $\norm{\xz} < 1$. Then the choice of $\delta = \bigl(\eps L^{-T}\bigr)\mn 1$ immediately gives us the (AS1) property.

			It remains to verify (SM1) and (SM2). Both the properties follow from~\eqref{e:supmartineq} in the proof of Proposition~\ref{p:supmartss}, with $s = 0$ and $\tau_r = \tau_{\{0\}} = 0$. Indeed, with these values of $s$ and $\tau_r$,~\eqref{e:supmartineq} becomes
			\begin{align*}
				\msf E\bigl[\epower{\alpha t}\alpha_1(\norm{X_t})\bigr|(\sigma_0, X_0)\bigr] & \le \msf E\bigl[\epower{\alpha t} V_{\sigma_t}(X_t)\bigr|(\sigma_0, X_0)\bigr] \le V_{\sigma_0}(X_0)(\alpha')^t \le \alpha_2(\norm{\xz})(\alpha')^t
			\end{align*}
			in view of (V1), where $\alpha' = \lambda_\circ(\hat p + \mu\tilde p)\epower\alpha < 1$. Therefore, given $\eps > 0$, we simply choose $\delta < \alpha_2^{-1}(\eps)$ to get (SM1). Given $r, \eps' > 0$, we simply choose $T = 0\mx \bigl(\ln(\alpha_2(r)/\eps')/\ln(\alpha')\bigr)$ to get (SM2). This completes the proof.
		\end{proof}

	\subsection{Robust Stability of Discrete-Time Randomly Switched Systems}\mbox{}
	\label{s:ifsd}
		Conditions for the existence of the supermartingale $\bigl(\epower{\alpha (t\mn\tauk)}V\bigl(X_{t\mn\tauk}\bigr)\bigr)_{t\in\Nz}$ in~\secref{s:gen} can be easily expressed in terms of the transition kernel $Q$. However, if $Q$ is not known exactly, which may happen if the model of the underlying system generating the Markov process $(X_t)_{t\in\Nz}$ is uncertain, one needs different methods. We look at one such instance below.

		Consider the system
		\begin{equation}
		\label{e:sysd}
			X_{t+1} = f_{\sigma_t}(X_t, w_t),\qquad X_0 = \xz, \quad t\in\Nz,
		\end{equation}
		where we retain the definition $\sigma$ from~\secref{s:ifs}, $f_i:\R^d\times\R^m\lra\R^d$ is locally Lipschitz continuous in both arguments with $f_i(0, 0) = 0$ for each $i\in\mcal P$, and $(w_t)_{t\in\Nz}$ is a bounded and measurable $\R^m$-valued disturbance sequence. We do not model $(w_t)_{t\in\Nz}$ as a random process; as such, the transition kernel of~\eqref{e:sysd} is not unique.

		\begin{defn}
		\label{d:issm}
			The system~\eqref{e:sysd} is said to be \emph{input-to-state stable in $\Lp 1$} if there exist functions $\chi, \chi'\in\ClassKinfty$ and $\psi\in\ClassKL$ such that $\EE_{\xz}\bigl[\chi(\norm{X_t})\bigr] \le \psi(\norm{\xz}, t) + \sup_{s\in\Nz}\chi'(\norm{w_s})$ for all $t\in\Nz$.\DefEnd
		\end{defn}

		%To wit, the system~\eqref{e:sysd} is input-to-state (\iss{}) stable in $\Lp 1$ if the $\Lp 1$ norm of the state is asymptotically bounded by a function of the size of the noise, and is uniformly bounded in expectation. 
		Our motivation for this definition comes from the concept of input-to-state stability \iss{} in the deterministic context~\cite{jiang01}. Consider the $i$-th subsystem of~\eqref{e:sysd} $x_{t+1} = f_i(x_t, w_t)$ for $t\in\Nz$ with initial condition $x_0$; note that $(x_t)_{t\in\Nz}$ is a deterministic sequence. This nonlinear discrete-time system is said to be \iss{} if there exist functions $\psi\in\ClassKL$ and $\chi\in\ClassKinfty$ such that $\norm{x_t} \le \psi(\norm{\xz}, t) + \sup_{s\in\Nz} \chi(\norm{w_s})$ for $t\in\Nz$. A sufficient set of conditions (cf.~\cite[Lemma~3.5]{jiang01}) for \iss{} of this system is that there exist a continuous function $V:\R^d\lra\posR$, $\alpha_1, \alpha_2\in\ClassKinfty$, $\rho\in\ClassK$, and a constant $\lambda\in\;]0, 1[$, such that $\alpha_1(\norm x) \le V(x) \le \alpha_2(\norm x)$ for all $x\in\R^d$, and $V(f_i(x, w)) \le \lambda V(x)$ whenever $\norm{x} > \rho(\norm w)$.

		In this framework we have the following Proposition.% observe that its hypotheses require that the subsystems of~\eqref{e:sysd} are sufficiently nice, in the sense that under zero-disturbance they reduce to Assumption~\ref{a:Vprime} with $r = 0$ (since $f_i(0, 0) = 0$).

		\begin{proposition}
			Consider the system~\eqref{e:sysd}, and suppose that
			\begin{enumerate}[label={\rm (\roman*)}, align=right, leftmargin=*, widest=iii]
				\item Assumption~\ref{a:sigma} holds,
				\item there exist continuous functions $V_i:\R^d\lra\posR$ for $i\in\PSet$, $\alpha_1, \alpha_2, \rho\in\ClassKinfty$, a constant $\mu > 1$ and a matrix $[\lambda_{ij}]_{\cardP\times\cardP}$ of nonnegative entries, such that
				\begin{enumerate}[label={\rm (\alph*)}, align=right, leftmargin=*]
					\item $\alpha_1(\norm{x}) \le V_i(x) \le \alpha_2(\norm x)\qquad$ for all $x$ and $i$,
					\item $V_i(x) \le \mu V_j(x)\qquad$ for all $x$ and $i, j$, and 
					\item $V_i(f_j(x)) \le \lambda_{ij} V_i(x)\qquad$ whenever $\norm{x} > \rho(\norm w)$ and all $i, j$,
				\end{enumerate}
				\item $\mu\Bigl(\max_{i\in\PSet}\sum_{j\in\PSet} p_{ij}\lambda_{ji}\Bigr) < 1$.
			\end{enumerate}
			Then~\eqref{e:sysd} is input-to-state stable in $\Lp 1$ in the sense of Definition~\ref{d:issm}.
		\end{proposition}
		\begin{proof}
			We define the compact set $K \Let  \bigl\{(i, y)\in\PSet\times\R^d\big| \norm{y} \le \sup_{s\in\Nz}\rho(\norm{w_s})\bigr\}$, and let $\tauk \Let  \inf\bigl\{t\in\Nz\big| X_t\in K\bigr\}$. In this setting we know from the preceding analysis that $\varphi(t, \xi) = \epower{\alpha t}\xi$, $\theta(t) = \epower{-\alpha t}$, and $C = 1/(1-\epower{-\alpha})$. We see from the estimate~\eqref{e:key} in the proof of Theorem~\ref{maintheorem} that 
			\begin{align*}
				\EE_{\xz}\bigl[V_{\sigma_t}(X_t)\bigr] & \le \varphi(0, V_{\sigma_0}(\xz))\theta(t) + \frac{\beta}{1-\epower{-\alpha}} + \delta \le \alpha_2(\norm{\xz})\epower{-\alpha t} + \frac{\beta}{1-\epower{-\alpha}} + \delta.
			\end{align*}
			Standard arguments show that there exists some $\chi''\in\ClassKinfty$ such that $\beta$ and $\delta$ are each dominated by $\chi''\bigl(\sup_{s\in\Nz}\norm{w_s}\bigr)$, and therefore, there exists some $\chi'\in\ClassKinfty$ such that $\beta/(1-\epower{-\alpha}) + \delta$ is dominated by $\chi'\bigl(\sup_{s\in\Nz}\norm{w_s}\bigr)$. Applying (ii)(a) on the left-hand side of the last inequality, we conclude that~\eqref{e:sysd} is input-to-state stable with $\chi = \alpha_1$ and $\psi(r, t) = \alpha_2(r)\epower{-\alpha t}$.
			%Let $\tau^1$ is the first entry time to the set $K$ of the process $(X_t)_{t\in\Nz}$ initialized at $\xz$; $\tau^1 > 0$ if $\xz\not\in K$, and equals $0$ otherwise. We observe that the Lyapunov functions satisfy the conditions in Assumption~\ref{a:Vprime}. As such, arguments identical to the proof of Proposition~\ref{p:supmartss} show that there exists $\alpha > 0$ such that the process $\bigl(\epower{\alpha(t\mn\tauk)} V_{\sigma_{t\mn\tauk}}\bigl(X_{t\mn\tauk}\bigr)\bigr)_{t\in\Nz}$ is a supermartingale whenever $\xz\not\in K$. Since $f_i$ is locally Lipschitz and $w$ is bounded, arguments similar to those in the proof of Corollary~\ref{c:gasas} show that $\beta \Let  \sup_{\xz\in K}\msf E\bigl[V_{\sigma_1}(X_1)\indic{\{X_1\in\R^d\setminus K\}}\big|X_0 = \xz\bigr] < \infty$. The proof of Theorem~\ref{maintheorem} now shows that $\msf E\bigl[V_{\sigma_t}(X_t)\indic{\{t\ge \tau^1\}}\bigr] \le \sup_{i\in\PSet, y\in K} V_i(y) + \beta \le \alpha_2\circ\rho\big(\sup_{s\in\Nz}\norm{w_s}\bigr) + \beta$. Standard arguments show that $\beta$ can be dominated by $\gamma\bigl(\sup_{s\in\Nz}\norm{w_s}\bigr)$. Also, if $\xz\in\R^d\setminus K$, then $\Ex{V_{\sigma_t}(X_t)\indic{\{t < {\tau}^1\}}} \le \alpha_2(\norm{\xz})\epower{\alpha t}$.
		\end{proof}

	\section*{Acknowledgments}
		The authors thank Daniel Liberzon and John Lygeros for helpful comments, Andreas Milias-Argeitis for useful discussions related to the chemical master equation, and the anonymous reviewer for a thorough review of the manuscript, several helpful comments, and drawing their attention to \cite[Chapter 14]{meynbk}.

\bibliographystyle{amsalpha}

%\bibliography{../../references}

\begin{thebibliography}{50}

%\bibitem[AS92]{aldenRH}
%J.~M. Alden and R.~L. Smith, \emph{Rolling horizon procedures in nonhomogeneous
%  {M}arkov decision processes}, Operations Research \textbf{40} (1992),
%  no.~suppl.~2, S183--S194.

\bibitem[ACK08]{kurtz08}
D.~F. Anderson, G. Craciun, and T.~G. Kurtz, \emph{Product-form stationary distributions for deficiency zero chemical reaction networks}, \url{http://arxiv.org/abs/0803.3042}, 2008.

%\bibitem[ABFG{\etalchar{+}}93]{borkarsurvey93}
%A.~Arapostathis, V.~S. Borkar, E.~Fern{\'a}ndez-Gaucherand, M.~K. Ghosh, and
%  S.~I. Marcus, \emph{Discrete-time controlled {M}arkov processes with average
%  cost criterion: a survey}, SIAM Journal on Control and Optimization
%  \textbf{31} (1993), no.~2, 282--344.

\bibitem[BDEG88]{barnsley88}
M.~F. Barnsley, S.~G. Demko, J.~H. Elton, and J.~S. Geronimo, \emph{Invariant
  measures for {M}arkov processes arising from iterated function systems with
  place-dependent probabilities}, Annales de l'Institut Henri Poincar\'e.
  Probabilit\'es et Statistique \textbf{24} (1988), no.~3, 367--394, Erratum in
  ibid., {\bf 24} (1989), no.\ 4, 589--590.

\bibitem[BS78]{bertsekasshreve78}
D.~P. Bertsekas and S.~E. Shreve, \emph{Stochastic Optimal Control: the
  Discrete-Time Case}, Mathematics in Science and Engineering, vol. 139,
  Academic Press Inc. [Harcourt Brace Jovanovich Publishers], New York, 1978.

\bibitem[Bor91]{borkarTopicsControlledMC}
V.~S. Borkar, \emph{Topics in Controlled {M}arkov Chains}, Pitman Research
  Notes in Mathematics Series, vol. 240, Longman Scientific \& Technical,
  Harlow, 1991.

\bibitem[BKR{\etalchar{+}}01]{borodin01}
A.~Borodin, J.~Kleinberg, P.~Raghavan, M.~Sudan, and D.~P. Williamson,
  \emph{Adversarial queuing theory}, Journal of the ACM \textbf{48} (2001),
  no.~1, 13--38.

\bibitem[CCCL08]{cccl08}
D.~Chatterjee, E.~Cinquemani, G.~Chaloulos, and J.~Lygeros, \emph{Stochastic
  control up to a hitting time: optimality and rolling-horizon implementation},
  \url{http://arxiv.org/abs/0806.3008}, 2008.

\bibitem[CHL09]{chl09}
D.~Chatterjee, P.~Hokayem, and J.~Lygeros, \emph{Stochastic receding horizon
  control with bounded control inputs: a vector space approach},
  \url{http://arxiv.org/abs/0903.5444}, 2009.

\bibitem[CRS71]{CRS}
Y.~S.~Chow, H.~Robbins, and D.~Siegmund \emph{Great Expectations:
The Theory of Optimal Stopping},
Houghton Mifflin Company Boston, 1971. 

\bibitem[CFM05]{costaMJLS}
O.~L.~V. Costa, M.~D. Fragoso, and R.~P. Marques, \emph{Discrete-time {M}arkov
  Jump Linear Systems}, Probability and its Applications (New York),
  Springer-Verlag, London, 2005.

\bibitem[DF99]{diaconisifs}
P.~Diaconis and D.~Freedman, \emph{Iterated random functions}, SIAM Review
  \textbf{41} (1999), no.~1, 45--76 (electronic).

\bibitem[DFMS04]{ref:doucetal04}
R.~Douc, G.~Fort, E.~Moulines, and P.~Soulier, \emph{Practical drift conditions
  for subgeometric rates of convergence}, The Annals of Applied Probability
  \textbf{14} (2004), no.~3, 1353--1377.

\bibitem[FK04]{foss04}
S.~Foss and T.~Konstantopoulos, \emph{An overview of some stochastic stability
  methods}, Journal of Operations Research Society of Japan \textbf{47} (2004),
  no.~4, 275--303.

\bibitem[HLR96]{hastad96}
J.~H{\aa}stad, T.~Leighton, and B.~Rogoff, \emph{Analysis of backoff protocols
  for multiple access channels}, SIAM Journal on Computing \textbf{25} (1996),
  no.~4, 740--774.

%\bibitem[HLL90]{hernandez-lerma1990}
%O.~Hern{\'a}ndez-Lerma and J.~B. Lasserre, \emph{Error bounds for rolling
%  horizon policies in discrete-time {M}arkov control processes}, IEEE
%  Transactions on Automatic Control \textbf{35} (1990), no.~10, 1118--1124.
%
%\bibitem[HLL93]{hernandez-lerma1993}
%\bysame, \emph{Value iteration and rolling plans for {M}arkov control processes
%  with unbounded rewards}, Journal of Mathematical Analysis and Applications
%  \textbf{177} (1993), no.~1, 38--55.

\bibitem[HLL96]{hernandez-lerma1}
O.~Hern{\'a}ndez-Lerma and J.~B. Lasserre, \emph{Discrete-Time {M}arkov Control Processes: Basic Optimality
  Criteria}, Applications of Mathematics, vol.~30, Springer-Verlag, New York,
  1996.

\bibitem[HLL99]{hernandez-lerma2}
\bysame, \emph{Further Topics on Discrete-Time {M}arkov Control Processes}, Applications of Mathematics,
  vol.~42, Springer-Verlag, New York, 1999.

\bibitem[Hob98]{hobson98}
D.~G. Hobson, \emph{Volatility misspecification, option pricing and superreplication via coupling}, The Annals of Applied Probability \textbf{8} (1998) no. ~1, 193--205.

\bibitem[JH07]{huisinga07}
T.~Jahnke and W.~Huisinga, \emph{Solving the chemical master equation for
  monomolecular reaction systems analytically}, Journal of Mathematical Biology
  \textbf{54} (2007), no.~1, 1--26.

\bibitem[JT01]{tweedie01}
S.~F. Jarner and R.~L. Tweedie, \emph{Locally contracting iterated functions
  and stability of {M}arkov chains}, Journal of Applied Probability \textbf{38}
  (2001), no.~2, 494--507.

\bibitem[JW01]{jiang01}
Z-P. Jiang and Y.~Wang, \emph{Input-to-state stability for discrete-time
  nonlinear systems}, Automatica \textbf{37} (2001), no.~6, 857--869.

\bibitem[KS08]{karatshreve}
I. Karatzas and S. Shreve, \emph{Brownian Motion and Stochastic Calculus}, 2 ed., Graduate Texts in mathematics, Springer, 2008.  

%\bibitem[Kal02]{kallenberg02}
%O.~Kallenberg, \emph{Foundations of Modern Probability}, 2 ed., Probability and
%  Its Applications, Springer Verlag, 2002.

\bibitem[Kif86]{kifer86}
Y.~Kifer, \emph{Ergodic Theory of Random Transformations}, Progress in
  Probability and Statistics, vol.~10, Birkh\"auser Boston Inc., Boston, MA,
  1986.

%\bibitem[Kus71]{kushnerIntroStochControl}
%H.~Kushner, \emph{Introduction to Stochastic Control}, Holt, Rinehart and
%  Winston, Inc., New York, 1971.

\bibitem[Lib03]{liberzonbk}
D.~Liberzon, \emph{Switching in Systems and Control}, Systems \& Control:
  Foundations \& Applications, Birkh{\"a}user, Boston, 2003.

\bibitem[LM94]{lasota94}
A.~Lasota and M.~C. Mackey, \emph{Chaos, Fractals, and Noise}, 2 ed., Applied
  Mathematical Sciences, vol.~97, Springer-Verlag, New York, 1994.

\bibitem[LM02]{lasotamyjak02}
A.~Lasota and J.~Myjak, \emph{On a dimension of measures}, Polish Academy of
  Sciences. Bulletin. Mathematics \textbf{50} (2002), no.~2, 221--235.

\bibitem[LS04]{lasotaszarek04}
A.~Lasota and T.~Szarek, \emph{Dimension of measures invariant with respect to
  the {W}a\.zewska partial differential equation}, Journal of Differential
  Equations \textbf{196} (2004), no.~2, 448--465.

%\bibitem[LB07]{lorenziAMMS}
%L.~Lorenzi and M.~Bertoldi, \emph{Analytical Methods for {M}arkov Semigroups},
%  Pure and Applied Mathematics, vol. 283, Chapman \& Hall/CRC,
%  Boca Raton, FL, 2007.

\bibitem[Mac01]{maciejowskibk}
J.~M. Maciejowski, \emph{Predictive Control with Constraints}, Prentice Hall,
  2001.

\bibitem[Mey08]{meynCTCN}
S.~P. Meyn, \emph{Control Techniques for Complex Networks}, Cambridge
  University Press, Cambridge, 2008.

\bibitem[MT09]{meynbk}
S.~P. Meyn and R.~L. Tweedie, \emph{Markov Chains and Stochastic Stability}, 2nd Ed.,
  Cambridge University Press, Cambridge, UK, 2009.

\bibitem[MA08]{miliasforthcoming}
A.~Milias-Argeitis, \emph{Fast simulation of the chemical master equation with iterated function systems}, In preparation, 2008.

\bibitem[PR99]{pemantle99}
R.~Pemantle and J.~S. Rosenthal, \emph{Moment conditions for a sequence with negative drift to be uniformly bounded in {$L\sp r$}}, Stochastic Processes
  and their Applications \textbf{82} (1999), no.~1, 143--155.

\bibitem[PS06]{PeskirShiryaev}
G.~Peskir and A.~N. Shiryaev, \emph{Optimal {S}topping and {F}ree-{B}oundary
  {P}roblems}, Lectures in Mathematics ETH {Z\"u}rich, Birkh\"auser Verlag,
  Basel, 2006.

%\bibitem[Piu97]{piunovski97}
%A.~B. Piunovskiy, \emph{Optimal Control of Random Sequences in Problems with Constraints}, Mathematics and its Applications, vol. 410, Kluwer Academic
%  Publishers, Dordrecht, 1997, With a preface by V. B. Kolmanovskii and A. N.
%  Shiryaev.

\bibitem[RY99]{ry99}
D.~Revuz and M.~Yor, \emph{Continuous Martingales and Brownian Motion}, 3 ed.,
  Grundlehren der Mathematischen Wissenschaften, vol. 293, Springer-Verlag,
  Berlin, 1999.

\bibitem[Sza06]{szarek06}
T.~Szarek, \emph{Feller processes on nonlocally compact spaces}, The Annals of
  Probability \textbf{34} (2006), no.~5, 1849--1863.

\bibitem[vS00]{vanderschaftHbk}
A.~{van der Schaft} and H.~Schumacher, \emph{An Introduction to Hybrid Dynamical Systems}, Lecture Notes in Control and Information Sciences, vol.
  251, Springer-Verlag London Ltd., London, 2000.

\bibitem[Wer05]{werner05}
I.~Werner, \emph{Contractive {M}arkov systems}, Journal of the London
  Mathematical Society. Second Series \textbf{71} (2005), no.~1, 236--258.

\bibitem[Wil06]{wilkinson06}
D.~J. Wilkinson, \emph{Stochastic Modelling for Systems Biology}, Chapman \&
  Hall/CRC Mathematical and Computational Biology Series, Chapman \& Hall/CRC,
  Boca Raton, FL, 2006.

\end{thebibliography}

\bigskip\bigskip

\end{document}